\documentclass[a4paper]{article}
\usepackage{amsthm,amsmath,amssymb}
\usepackage{color}
\newtheorem{teor}{Theorem}[section]
\newtheorem{defin}[teor]{Definition}
\newtheorem{lemm}[teor]{Lemma}
\newtheorem{osse}[teor]{Remark}
\newtheorem{prop}[teor]{Proposition}
\newtheorem{defi}[teor]{Definition}
\newtheorem{coro}[teor]{Corollary}
\newtheorem{prob}[teor]{Problem}

\newcommand{\bele}{\begin{lemm}\begin{sl}}
\newcommand{\enle}{\end{sl}\end{lemm}}
\newcommand{\bedef}{\begin{defi}\begin{sl}}
\newcommand{\eddef}{\end{sl}\end{defi}}
\newcommand{\bete}{\begin{teor}\begin{sl}}
\newcommand{\ente}{\end{sl}\end{teor}}
\newcommand{\beos}{\begin{osse}\begin{rm}}
\newcommand{\eddos}{\end{rm}\end{osse}}
\newcommand{\bepr}{\begin{prop}\begin{sl}}
\newcommand{\empr}{\end{sl}\end{prop}}
\newcommand{\bepro}{\begin{prob}\begin{rm}}
\newcommand{\empro}{\end{rm}\end{prob}}
\newcommand{\bede}{\begin{defin}\begin{sl}}
\newcommand{\edde}{\end{sl}\end{defin}}
\newcommand{\beco}{\begin{coro}\begin{sl}}
\newcommand{\enco}{\end{sl}\end{coro}}

\definecolor{grey}{rgb}{0.85,0.85,0.85}

\long\def\greybox#1{%
    \newbox\contentbox%
    \newbox\bkgdbox%
    \setbox\contentbox\hbox to \hsize{%
        \vtop{
            \kern\columnsep
            \hbox to \hsize{%
                \kern\columnsep%
                \advance\hsize by -2\columnsep%
                \setlength{\textwidth}{\hsize}%
                \vbox{
                    \parskip=\baselineskip
                    \parindent=0bp
                    #1
                }%
                \kern\columnsep%
            }%
            \kern\columnsep%
        }%
    }%
    \setbox\bkgdbox\vbox{
        \color{grey}
        \hrule width  \wd\contentbox %
               height \ht\contentbox %
               depth  \dp\contentbox
        \color{black}
    }%
    \wd\bkgdbox=0bp%
    \vbox{\hbox to \hsize{\box\bkgdbox\box\contentbox}}%
    \vskip\baselineskip%
}

\newcommand{\RR}{\mathbb{R}}

\newcommand{\beeq}[1]{\begin{equation}\label{#1}}
\newcommand{\eddeq}{\end{equation}}

\newcommand{\beeqa}[1]{\begin{eqnarray}\label{#1}}
\newcommand{\eddeqa}{\end{eqnarray}}

\newcommand{\beal}[1]{\begin{align}\label{#1}}
\newcommand{\eddal}{\end{align}}

\newcommand{\bespl}[1]{\begin{split}\label{#1}}
\newcommand{\edspl}{\end{split}}

\newcommand{\bega}[1]{\begin{gather}\label{#1}}
\newcommand{\edga}{\end{gather}}

\newcommand{\beeqax}{\begin{eqnarray*}}
\newcommand{\eddeqax}{\end{eqnarray*}}

\def\qed{\ifmmode 
  \else \leavevmode\unskip\penalty9999 \hbox{}\nobreak\hfill
  \fi
  \quad\hbox{\hskip.5em\vrule width.4em height.6em depth.05em\hskip.1em}}
\def\endproofsym{\qed}
\renewenvironment{proof}[1][Proof]{\trivlist\item[\hskip\labelsep{\hskip0pt
    {\normalfont\scshape#1.}\hskip .321429\parindent}]\ignorespaces}
{\endproofsym\endtrivlist}
\def\endnobox{\def\endproofsym{}\end{proof}\def\endproofsym{\qed}}

\newcommand{\no}{\nonumber}

\newcommand{\beeqao}{\begin{eqnarray}\no}
\newcommand{\bealo}{\begin{align}\no}
\newcommand{\besplo}{\begin{split}\no}
\newcommand{\begao}{\begin{gather}\no}

\newcommand{\ep}{\varepsilon}
\newcommand{\eps}{\varepsilon}

\newcommand{\ov}[1]{\overline{#1}}
\newcommand{\itt}{\int_0^t}

\newcommand{\io}{\int_\Omega}

\newcommand{\ito}{\itt\!\io}

\newcommand{\iTT}{\int_0^T}

\newcommand{\iTo}{\iTT\!\io}

\newcommand{\ee}{_{\varepsilon}}

\newcommand{\lhs}{left hand side}
\newcommand{\rhs}{right hand side}

\newcommand{\TS}{R^{3 \times 3}_{{\rm sym},0} }

\DeclareMathOperator{\tr}{tr}

\DeclareMathOperator{\deriv}{d}

\DeclareMathOperator{\sign}{sign}

\DeclareMathOperator{\ess}{ess}

\let\TeXchi\chi
\def\chi{{\setbox0 \hbox{\mathsurround0pt
$\TeXchi$}\hbox{\raise\dp0 \copy0 }}}

\newcommand{\teta}{\vartheta}

\newcommand{\kappaciapo}{\widehat{\kappa}}

\newcommand{\calF}{{\mathcal F}}

\newcommand{\mugiu}{\underline{\mu}}
\newcommand{\musu}{\overline{\mu}}

\newcommand{\dit}{\deriv\!t}

\newcommand{\ddt}{\frac{\deriv\!{}}{\dit}}

\newcommand{\calL}{{\cal L}}

\newcommand{\bFormula}[1]{
\begin{equation} \label{#1}}
\newcommand{\eF}{\end{equation}}

\newcommand{\DC}{C^\infty_c}

\newcommand{\vt}{\vartheta}
\newcommand{\vu}{\vc{u}}

\newcommand{\vc}[1]{{\bf #1}}
\newcommand{\Div}{{\rm div}_x}
\newcommand{\Grad}{\nabla_x}

\newcommand{\HH}{\mathbb{H}}
\newcommand{\II}{\mathbb{I}}
\newcommand{\QQ}{\mathbb{Q}}
\newcommand{\VV}{\mathbb{V}}
\newcommand{\MM}{\mathbb{M}}
\newcommand{\SSS}{\mathbb{S}}

\newcommand{\pp}{\mathbf{p}}
\newcommand{\calD}{{\mathcal D}}

\def\bproof{\noindent{\bf Proof.\;}}
\def\be{\begin{equation}}
\def\ee{\end{equation}}
\def\eproof{\hfill$\square$\medskip}

\def\non{{\nonumber}}
\def\Sphere{{\mathbb S}}

\def\bZ{{\bar Z}}

\def\mcI{{\mathcal{I}}}

\newcommand{\emu}{e^{3\mu\sin^2\theta}}

\numberwithin{equation}{section}
\begin{document}

\title{Nonisothermal nematic liquid crystal flows with
 the Ball-Majumdar free energy}

\author{%
Eduard Feireisl\\
Institute of Mathematics of the Czech Academy of Sciences,\\
\v Zitn\' a 25, 115 67 Praha 1, Czech Republic\\
E-mail: {\tt feireisl@math.cas.cz}
\thanks{The research of E.F. leading to these results has received funding from the European Research Council under the European Union's Seventh Framework
Programme (FP7/2007-2013)/ ERC Grant Agreement 320078.The work of E.F. was partially supported by
  the FP7-IDEAS-ERC-StG \#256872 (EntroPhase).}
\and
Giulio Schimperna\\
Dipartimento di Matematica, Universit\`a di Pavia,\\
Via Ferrata~1, 27100 Pavia, Italy\\
E-mail: {\tt giusch04@unipv.it}%
\thanks{The work of G.S. was supported by the
 FP7-IDEAS-ERC-StG \#256872 (EntroPhase).} \\
\and
Elisabetta Rocca\\
Weierstrass Institute for Applied
Analysis and Stochastics, \\
Mohrenstr.~39, D-10117 Berlin, Germany\\
E-mail: {\tt  rocca@wias-berlin.de}\\
and \\
Dipartimento di Matematica, Universit\`a di Milano,\\
Via Sal\-di\-ni 50, 20133 Milano, Italy\\
E-mail {\tt elisabetta.rocca@unimi.it}%
\thanks{The work of E.R. was supported by the
 FP7-IDEAS-ERC-StG \#256872 (EntroPhase).}
\and
Arghir Zarnescu\\
Pevensey III, University of Sussex,\\
Falmer, BN1 9QH, UK\\
E-mail: {\tt A.Zarnescu@sussex.ac.uk}%
\thanks{The work of A.Z. was partially supported by
  the FP7-IDEAS-ERC-StG \#256872 (EntroPhase).}
}

\date{}

\maketitle
\begin{abstract}
In this paper we prove  the existence of global in time weak solutions for an evolutionary
PDE system modelling nonisothermal Landau-de Gennes nematic liquid crystal (LC) flows
in three dimensions of space. In our model, the incompressible Navier-Stokes system
for the macroscopic velocity $\vu$ is coupled to a nonlinear convective parabolic equation
describing the evolution of the Q-tensor $\QQ$, namely a tensor-valued variable
representing the normalized second order moments of the probability
distribution function of the  LC molecules. The effects of the (absolute)
temperature $\vt$ are prescribed in the form of an energy balance
identity complemented with a global entropy production inequality.
Compared to previous contributions, we can consider here the
physically realistic singular configuration potential $f$ introduced by
Ball and Majumdar. This potential
gives rise to severe mathematical difficulties since it introduces,
in the Q-tensor equation, a term which is at the same time singular
in $\QQ$ and degenerate in $\vt$. To treat it a careful
analysis of the properties of $f$, particularly of its blow-up
rate, is carried out.
\end{abstract}


\vspace{2mm}

\noindent {\bf AMS (MOS) subject clas\-si\-fi\-ca\-tion:}%
~~76A15, 74G25, 35D30, 35Q30.


\section{Introduction}
\label{sec:intro}

The main aim of this paper is to continue the analysis of non-isothermal Landau-de Gennes nematic
liquid crystal flows with singular potential started in our recent contribution~\cite{FRSZ}.
Our purpose is to consider the evolution of the flow assuming that the mixing term in the
bulk potential is singular, and, in particular, takes the form introduced by Ball and Majumdar
in \cite{BalMaj}:
\bFormula{psimajIntro}
  \psi_B (\vt, \QQ) = \vt f(\QQ)+ G(\QQ),
\eF
where $\vt$ denotes the absolute temperature and the {\sl Q-tensor}\/
$\QQ$ is a symmetric traceless tensor representing in a suitable way the
normalized second order moments of the probability
distribution function of molecules. Here $G$ is a  smooth function of $\QQ$
while $f: R^{3 \times 3}_{\rm sym,0} \to [-K,+\infty]$,
$K\ge 0$, is a convex function smooth on its domain which consists of $\QQ$-tensors whose eigenvalues are  inside the {\sl physical interval}~$(-1/3,2/3)$.

The use of Q-tensors permits to capture fine properties of the crystal configuration (
like biaxiality), which  are not possible to describe in the framework
of vector-based models. In order to introduce Q-tensors mathematically,
we start from a probability measure $\mu_x$ on the unit sphere $\mathbb{S}^2$,
representing the orientation of the molecules at a point $x$ in space. Then, we can
associate to $\mu_x$ a symmetric and traceless $3\times 3$ matrix defined as
\begin{equation}\nonumber
  \QQ(x)=\int_{\mathbb{S}^2}\left(\pp\otimes \pp-\frac{1}{3}\II\right)\ {\rm d}\mu_x(\pp).
\end{equation}
Namely, $\QQ$ is intended to measure how much the probability measure $\mu_x$ deviates
from the isotropic measure $\bar\mu$ where $d\bar\mu=\frac{1}{4\pi}dA$, see \cite{DeGennes}.
In the Onsager model (cf.~\cite{DeGennes}, \cite{MaSa}), $\mu_x$ is
assumed to be absolutely continuous with density $\rho=\rho(\pp)$. In
this case, we have
\begin{equation}\label{Qrep}
  \QQ(x)=\int_{\mathbb{S}^2}\left(\pp\otimes \pp-\frac{1}{3}\II\right)\,\rho(\pp)
  \ {\rm d}\pp.
\end{equation}
The fact that $\mu_x$ is a probability measure imposes a constraint on the eigenvalues of
$\QQ$, namely that they are bounded between the values $-1/3$ and $2/3$,
see \cite{BalMaj}. Thus not any traceless $3\times 3$ matrix is a {\it physical} $\QQ$-tensor
but only those whose eigenvalues are in $(-1/3,2/3)$.

The Ball-Majumdar potential is then defined as follows:
\bFormula{f1}
  \displaystyle f(\QQ) = \left\{ \begin{array}{l} \inf_{\rho \in \mathcal{A}_{\QQ}}
   \int_{S^2} \rho(\vc{p}) \log (\rho(\vc{p})) \ {\rm d}\vc{p}
    \ \mbox{if} \ \lambda_i[\QQ] \in (-1/3,2/3),\ i=1,2,3,\\
    \\ + \infty \ \mbox{otherwise,}
         \end{array} \right.
\eF
\[
 \displaystyle  \mathcal{A}_{\QQ}
    = \left\{ \rho: S^2 \to [0, \infty) \
     \Big|\ \rho\in L^1(S^2),\ \int_{S^2} \rho (\vc{p}) \ {\rm d} \vc{p} = 1;
    \QQ = \int_{S^2} \left( \vc{p} \otimes \vc{p} - \frac{1}{3} \II \right) \rho(\vc{p}) \ {\rm d} \vc{p}
        \right\}.
\]
Then, it turns out (see, e.g., \cite{BalMaj})
that the {\sl effective domain}\/
$\calD[f]$ (i.e., the set where $f$ assumes finite values) coincides precisely
with the set of symmetric traceless tensors whose eigenvalues belong to
the physical interval. Moreover, as we prove in detail in
Section~\ref{sec:BM} below, $f$ explodes logarithmically fast as
$\QQ$ approaches the boundary of $\calD[f]$ (or, equivalently, one
of its eigenvalues  tends to $-1/3$). The main mathematical features
of $f$, which are required in the existence proof for the
purpose of obtaining the necessary a-priori estimates,
are thoroughly discussed in Section~\ref{sec:BM}.
Here, it is just worth noting that,
for large values of the temperature, the convex part $f$ of the energy
$\psi_B$ is prevailing, whereas, when $\vt$ is close to $0$,
$\psi_B$ may exhibit a multiple well structure.

Dealing with the  bulk potential \eqref{psimajIntro} generates severe
mathematical difficulties, since the term $\vt f(\QQ)$ is both
singular in $\QQ$ (since the eigenvalues of $\QQ$ may approach,
at least in some small region, the limiting values $-1/3,2/3$,
and degenerate in $\vt$ (since we can only prove that $\vt>0$
almost everywhere, not excluding, however, that  the (essential) infimum 
of $\vt$
might reach $0$). Because of several technical problems associated with the above mentioned difficulties, our previous contribution \cite{FRSZ} 
focused on
the case when the bulk potential
in the free energy functional is given by (cf. also \cite{SchoSlu} and
\cite{MottramNewton})
\bFormula{psi}
  \psi_B (\vt, \QQ) = f(\QQ) - U(\vt) G(\QQ),
\eF
where $U$ is assumed to be a convex and decreasing function of $\vt$ with
restricted growth at $\infty$. This form of the bulk potential allowed us to get rid of the degenerate character
in $\vt$, and,  in fact, to prove strict positivity of $\vt$
in that case. Although the formula \eqref{psimajIntro} is \emph{mathematically} equivalent to $\eqref{psi}$ in the framework of stationary 
problems, they give rise to different solutions if the time evolution of the system is studied. In accordance with the underlying physical arguments 
(see \cite{MaSa}), it is \eqref{psimajIntro} rather than \eqref{psi} that should be considered. 
Indeed this choice, combined with the standard principles of Thermodynamics, yields 
the entropy of the system in the form
\[
  s = - \frac{\partial \mathcal{F}}{\partial \vt}
   = c(\vt)- f(\QQ),
\]
where $\mathcal{F}$ denotes the free energy functional, and $c$ represents
its purely caloric part (cf.~also \eqref{psi}). On the other hand,
in  the case \eqref{psi}, one has
\[
  s=c(\vt)+U'(\vt)G(\QQ),
\]
and the singular potential $f$ appears only in the internal energy of the system.

In our model, the evolution of the Q-tensor is derived as a balance of
the following free energy functional:
\bFormula{i1}
  \mathcal{F} (\QQ, \Grad \QQ, \vt )
   = \frac{1}{2} |\Grad \QQ |^2
    + \psi_B (\vt, \QQ)
    - \vt \log \vt
    - a \vt^m,
\eF
where $\lambda,a\ge 0$ and the exponent $m\ge 1$ will be fixed below.
The term depending on $\Grad \QQ$  describes the interfacial energy,
whereas $a \vt^m$, $a>0$, prescribes a power-like heat conductivity.

The macroscopic velocity of the crystal flow is described by
a Navier-Stokes type system (cf.~\eqref{i22} below) for the velocity field $\vu$, in
which the stress tensor contains a non-standard part depending on the tensor~$\QQ$.
The evolution of $\QQ$ (cf.~\eqref{i23}) is ruled by the variational derivative
of the free-energy functional, displaying a singular character due to the presence
of the singular Ball-Majumdar energy \eqref{psimajIntro}.
In the physical
formulation of the model, the evolution of temperature is described by
the heat equation \eqref{p6}. However, due to the presence of quadratic
terms in the \rhs\ (especially those depending on $\vu$), dealing
with \eqref{p6} appears out of reach from a mathematical perspective.
For this reason, following an
idea originally developed by Bul\'\i\v cek, Feireisl and M\'alek
in \cite{BFM} for the non-isothermal Navier-Stokes system, we
replace it by an equality describing conservation of total energy conservation,
combined with a weak form of the entropy inequality (cf.~\eqref{i24}
and \eqref{i28weak} below). As a drawback of this choice,
we have to note the explicit appearence of the pressure $p$ in
the energy balance equation. Actually, $p$ is no longer just a Lagrange
multiplier, but it enters the system as an additional unknown
and needs to be controlled carefully. In particular, for this purpose
a suitable choice of the boundary condition is crucial and,
indeed, in order to reduce complications to the minimum, we consider
periodic b.c.'s (cf.~\eqref{toro} below).

Due to the quoted mathematical difficulties, the results
we obtain in the present contribution are {\sl weaker}\/ than
those we got in \cite{FRSZ}. Actually, here we consider
a weaker notion of solution, where the entropy inequality
holds only in its integrated (both in space ad in time)
form; in other words, we can only prove a {\sl global balance}\/ of entropy.
In addition, in order to deal with the Ball-Majumdar energy
\eqref{psimajIntro}, we need a power-like growth assumptions imposed on
the heat conductivity, the specific heat and the collective rotational
viscosity coefficient in the system. Mathematically
speaking, these assumptions are enforced by the need of sufficiently strong
{\it a-priori} bounds for the temperature $\vt$.
The derivation of suitable {\it a priori} estimates as well as
a rigorous justification of compactness of the family of solutions
represents the main novelty of the present paper with respect to \cite{FRSZ}, while the 
specific construction of suitable approximate solutions, carried over in detail in \cite{FRSZ}, 
is a routine matter. This is the reason why we restrict our discussion to the passage  from {\it a priori} bounds to  compactness 
(sequential stability), leaving the necessary modifications of the construction of suitable approximate solutions in \cite{FRSZ} to the interested
reader.

As an additional result, in the last part of the paper
we see that if we additionally assume that the heat flux is {\sl singular}
as the absolute temperature $\vt$ approaches $0$
(see Section \ref{sec:singheat} below), then we can recover the
entropy inequality in the usual distributional sense. Namely, we can
estimate the entropy growth on any region of the space domain,
and not only globally. Moreover, we can deal with an even more
general class of singular potentials. It is worth noting that
taking this type of singular flux law is not only a mathematically
ad-hoc assumption, but is a common choice in several types of
phase-transition and phase-separation models both in liquids
and in crystalline solids (cf., e.g., \cite{colli-lau}, \cite{RS}, and \cite{PF}
where similar growth conditions are assumed).
%
%
Regarding the choice of power-growing specific heat (cf.~the last summand
in \eqref{i1}), let us mention the papers \cite{AP, MS} for examples in
phase transitions and \cite{BG, rocca-rossi-deg, rossi-roubi} for examples in thermoelastic systems, where
this behaviour is observed. Note that, however, in our contribution we can allow the viscosity coefficient to depend on the absolute temperature in a quite general way: it has to be only bounded from below and above without any specific growth condition. 

We conclude this introduction by a brief review of other  LC models
based either on the tensorial variable $\QQ$ or other parameters. To begin, a vectorial quantity ${\bf d}$ could be used
instead of $\QQ$ in order to describe the preferred orientation of the molecules
at any point (cf., e.g., \cite{CR, E91, LL01} and references therein).  Moreover, we can quote the standard
hydrodynamic theories, where the case of regular bulk potential
$$
  \psi_B(\QQ)=\frac{a}{2}{\rm tr}(\QQ^2)-\frac{b}{2}{\rm tr}(\QQ^3)+\frac{c}{4}{\rm tr}^2(\QQ^2)
$$
is considered (cf.~\cite{PaZa2011} and \cite{PaZa2012} for the analysis of the corresponding
isothermal evolution),  it is not clear whether this {\sl physical}\/ contraint on the
eigenvectors is preserved or not, while the choice of the Ball-Majumdar free energy
\eqref{psimajIntro} naturally enforces this constraint.
Let us notice that in the literature there are few papers dealing with liquid crystal
models with singular potential. In the isothermal case, we can quote the very recent
contribution \cite{Wilk}, where existence and regularity of solution in the 3D and 2D
cases are obtained for a tensorial model with Ball-Majumdar potential. Regarding the
non-isothermal tensorial case, up to our knowledge, the only contribution, in the case
of singular potential \eqref{psi}, was given in the paper \cite{FRSZ},
which we already mentioned in this Introduction. Two attempts of considering the
non-isothermal case were made in~\cite{FRS} and~\cite{FFRS} and refer to
vectorial models. In particular, in~\cite{FFRS} the stretching and rotation
effects of the director field induced by the straining of the fluid were
considered and the existence of global in time weak
solutions was obtained for the corresponding initial boundary value problem.

\bigskip

\noindent
{\bf Plan of the paper.}~~%
The remainder of the paper is organized as follows:
in Section~\ref{sec:model}
we first introduce the equations of our mathematical
model in the strong (physical) form (Subsection~\ref{subsec:eqn}); then,
we specify our basic assumptions on coefficients and data
(Subsection~\ref{subsec:hyp}). This permits us to
introduce (Subsection~\ref{subsec:weak}) a weak formulation
of the initial-boundary value problem for our model and
to present a rigorous existence theorem (Theorem~\ref{th:main})
 of weak solutions. The proof is carried out in Section~\ref{sec:proof}
and is subdivided into several steps: we first prove (Subsection~\ref{sec:apriori})
that any weak solution complies, at least formally, with a number of a-priori
estimates; subsequently, these bounds are shown to be sufficient
to provide weak sequential stability (Subsection~\ref{sec:limit}).
Namely, sequences of solutions starting from (suitably) bounded
families of initial data admit weak limits which still solve the
problem. A key ingredient in the proof is represented by some
fine properties of the Ball-Majumdar potential which are  obtained separately
in Section~\ref{sec:BM}. Finally, in Section~\ref{sec:singheat}
we show that, assuming a singular heat flux law, even more general
classes of potentials can be treated and, moreover,
the entropy inequality holds in a stronger sense.


\section{The mathematical problem}
\label{sec:model}


\subsection{The equations of the model}
\label{subsec:eqn}

In this part, we write the equations of our mathematical model
describing the evolution of the unknown fields $\vu$ (macroscopic
velocity), $\QQ$ (Q-tensor) and $\vt$ (absolute temperature).
The model can be physically derived by following closely the
argument reported in \cite[Section~1]{FRSZ}.
Indeed, the main  difference is represented here by the
choice of the Ball-Majumdar potential \eqref{psimajIntro}.
However, this does not affect the way the model is obtained,
but only the outcoming form of the equations.
These are stated here in their strong form.
In Section~\ref{subsec:weak}, we will see that, in order to have
a mathematically tractable system, passing to a
weak formulation is actually necessary.

In order to avoid complications related to interactions with
boundaries, we take periodic boundary conditions for all
unknowns. Namely, we assume the system be settled in the
three dimensional torus
\bFormula{toro}
  \Omega := \left( [-\pi, \pi]|_{\{ -\pi, \pi \} } \right)^3.
\eF
As already noted  in the Introduction, the above choice is crucial for our existence
theorem.

The first two relations in our model describe the evolution
of the macroscopic velocity $\vu$:

\medskip

\greybox{
\centerline{\textsc{Incompressibility:}}

\vspace{-6mm}

\bFormula{i21}
  \Div \vu = 0.
\eF
}

\smallskip

\greybox{
\centerline{\textsc{Momentum equation:}}

\vspace{-6mm}

\bFormula{i22}
  \partial_t \vu + \Div (\vu \otimes \vu)
   = \Div \sigma.
\eF
}
\noindent
Here, $\sigma$ denotes the {\sl stress tensor}, given by
\bFormula{i25}
  \sigma = \mu(\vt) \left( \Grad \vu + \Grad^t \vu \right) - p \II
\eF
\[
  \mbox{} + 2 \xi \left( \HH: \QQ \right) \left( \QQ + \frac{1}{3} \II \right)
   - \xi \left[ \HH \left( \QQ  + \frac{1}{3} \II \right)
   + \left( \QQ  + \frac{1}{3} \II \right) \HH \right]
   + \left( \QQ \HH - \HH \QQ \right)
   - \Grad \QQ \odot \Grad \QQ,
\]
$\xi$ being a fixed scalar parameter, measuring the ratio between the rotation
and the aligning effect that a shear flow exerts over the directors, see
Beris and Edwards \cite{BeEd}.

The behavior of the Q-tensor is ruled by the following relation:

\medskip

\greybox{
\centerline{\textsc{Order parameter evolution:}}

\vspace{-6mm}

\bFormula{i23}
  \partial_t \QQ
   + \vu \cdot \Grad \QQ
   - \SSS (\Grad \vu, \QQ)
  = \Gamma(\vt) \HH.
\eF
}
\noindent
The nonnegative function $\Gamma(\vt)$ represents a  collective rotational viscosity coefficient,
whereas
\bFormula{i5}
  \SSS(\Grad \vu, \QQ)
   = \left( \xi \ep(\vu) + \omega (\vu) \right) \left( \QQ + \frac{1}{3} \II \right)
    + \left( \QQ + \frac{1}{3} \II \right) \left( \xi \ep (\vu) - \omega (\vu) \right)
    - 2 \xi \left( \QQ + \frac{1}{3} \II \right) \left( \QQ: \Grad \vu \right)
\eF
denotes stretching and rotation effects on the Q-tensor driven by the
macroscopic flow, with $\omega(\vu)$ standing for the
antisymmetric part of the viscous stress, i.e.,
\[
  \omega(\vu) = \frac{1}{2} \left( \Grad \vu - \Grad^t \vu \right),
\]
 and $\ep(\vu)$ for its symmetric part, i.e., 
\[
\ep(\vu)=\frac{1}{2}\left(\Grad\vu+\Grad^t\vu\right)\,. 
\]
Finally, $-\HH$ indicates the first variation of the free energy $\calF$
with respect to $\QQ$, namely
\bFormula{i20}
  \HH
   = \Delta \QQ
   - \vt \mathcal{L} \left[ \frac{\partial f(\QQ)}{\partial \QQ } \right]
   + \lambda \QQ, \ \ \lambda\ge 0,
\eF
 where we have chosen $G(Q)=-\lambda Q$ in \eqref{psimajIntro}.
The operator $\calL[ \cdot ]$ denotes projection on the space of
traceless tensor, that is $\calL [ \MM ] = \MM - \frac13 \tr(\MM)$
for $\MM\in R^{3\times 3}$. In other words, the expression
$\mathcal{L} \left[ \frac{\partial f(\QQ)}{\partial \QQ } \right]$
represents the subdifferential of $f$, with respect to the
structure of the space of traceless tensors, evaluated at the
point $\QQ$. Finally, the evolution of temperature is represented
{\sl in the strong formulation}\/ by the

\medskip

\greybox{
\centerline{\textsc{Heat equation:}}

\vspace{-6mm}

\bFormula{p6}
  \partial_t \vt  + a (m -1) \partial_t \vt^m
   + \vu \cdot \Grad \vt + a (m -1) \vu \cdot \Grad \vt^m
   - \Div \big( \kappa(\vt) \Grad \vt \big)
\eF
\[
  = \vt \big( \partial_t f (\QQ)
   + \vu \cdot \Grad f (\QQ ) \big)
   + \frac{\mu(\vt)}{2} \big| \Grad \vu + \Grad^t \vu \big|^2
   + \Gamma(\vt) | \HH |^2.
\]
}
\noindent%
The model can be derived from the free energy functional
\eqref{i1} by stating the balances of total energy and of
entropy. The procedure is completely analogous to that performed
in \cite[Section~1.3]{FRSZ},
to which we refer the reader for details.
We just note here that \eqref{p6} corresponds
to assuming the internal energy flux
being given by (cf.~\cite[formula~(1.30)]{FRSZ})
\bFormula{i26}
  \vc{q} = - \kappa (\vt) \Grad \vt -
    \Grad \QQ : \SSS(\Grad \vu, \QQ).
\eF
%


\subsection{Assumptions on coefficients and data}
\label{subsec:hyp}

We state here our basic assumptions on coefficients
and nonlinear terms in system
\eqref{i21}-\eqref{i26}, as well as on the initial
data. First of all, we ask the viscosity $\mu$ to
be bounded with respect to $\vt$, while we prescribe a power-like
law for the heat flux and a linear growth
of the diffusion coefficient $\Gamma$ at infinity:
\begin{align}\label{h11}
  & \mu \in C^1(\RR;\RR), \ \ 0 < \mugiu \le \mu(r) \le \musu \ \
   \text{for all }r\in\RR,\\
 \label{h12}
  & \kappa(r) = A_0 + A_k r^k, \ \ A_0, A_k>0, \\[1mm]
 \label{h13}
  & \Gamma(r) = \Gamma_0 + \Gamma_1 r, \ \ \Gamma_0, \Gamma_1>0.
\end{align}
In addition to these, we need a power-like heat conductivity
in equation \eqref{p6}. Namely, we take $a>0$ and
assume that the exponents $k$ in \eqref{h12}
and $m$ in \eqref{p6} are strictly positive and satisfy
the relation
\begin{equation}\label{defiA}
  A:=\frac{3k+2m}3 > 9, \quad \frac32<m\leq \frac{6k}{5}.
\end{equation}
Next, we let the initial data $\vu_0$, $\QQ_0$, and $\vt_0$
satisfy
\begin{align}\label{m6}
  & \vu_0 \in L^2(\Omega; R^3), \ \Div \vu_0 = 0, \\
 \label{m6b}
  & \QQ_0 \in H^1(\Omega; \TS), \ f(\QQ_0) \in L^1(\Omega),\\
 \label{m6c}
  & \vt_0 \in L^\infty(\Omega),\ \ess \inf_\Omega \vt_0 =: \underline{\vt} > 0.
\end{align}
Finally, we assume that $f$ satisfies at least the following
basic properties:
\begin{itemize}
 \item[(i)] $f : \TS \to [-K, +\infty]$ is convex and lower semi-continuous,
 with $K \ge 0$.
 \item[(ii)] The domain of $f$,
 \[
   \mathcal{D}[f] = \big\{ \QQ \in R^{3 \times 3}_{\rm sym,0} \ | \ f(\QQ) < + \infty \big\}
   = \big\{  \QQ \in R^{3 \times 3}_{\rm sym,0} \ | \ \lambda_i[\QQ] \in (- 1/3,2/3) \big\},
 \]
 is an open, convex, and bounded subset of $\TS$.
 \item[(iii)] $f$ is smooth in $\mathcal{D}[f]$, $\lambda\geq0$.
\end{itemize}
We will see in Section~\ref{sec:BM} below (cf.~also
\cite[Section 3]{BalMaj}) that (i)-(iii) hold in the case
when $f$ is given by the Ball-Majumdar potential
introduced in \eqref{f1}. This is a nontrivial issue since
relation \eqref{f1} describes $f$ implicitly in terms of the
orientation distribution(s) $\rho(\vc{p})$ generating $\QQ$.
Hence, the properties of $f$  as a function of $\QQ$
need to be properly demonstrated. On the other hand,
(i)-(iii) do not characterize the Ball-Majumdar potential
completely and are satisfied by a much wider class of convex
singular functions. Indeed, in the proof we will need
that $f$ satisfies some finer properties, especially related
to its blow  up rate as $\QQ$ approaches the boundary of
$\calD[f]$. The discussion of these properties,
and the proof of the fact that they are satisfied by
the Ball-Majumdar  singular potential, are postponed to
Section~\ref{sec:BM}.


\subsection{Weak formulation  and main results}
\label{subsec:weak}

A proof of existence of global-in-time \emph{strong} solutions
 of  the model introduced in
Section~\ref{subsec:eqn} appears to be out of reach for the available mathematical tools, one of
the principal difficulties
being related to the presence of quadratic terms on the
\rhs\ of the heat equation \eqref{p6}. Accordingly, 
we introduce here a weak formulation of the system
along the lines of \cite{FRSZ}. The main point stands
in the replacement of \eqref{p6} by means of a
couple of relations, namely the {\sl total energy balance}\/
(cf.~\eqref{i24} below) and a weak form of the
{\sl entropy production inequality}\/ (cf.~\eqref{i28weak}).

Such relations account for most of the information
coming from the heat equation, but are more treatable
mathematically. Actually, \eqref{i24} contains
no quadratic terms, while inequality \eqref{i28weak}
just states that the {\sl global}\/ entropy of the
system needs to grow (at least) as prescribed by the \rhs.
As noted in the Introduction, differently from
what was done in \cite{FRSZ} and due to technical difficulties
related to the energy~\eqref{psimajIntro}, we are not able
to get here any ``pointwise'' control of the entropy
(actually, this can be reached under different
assumptions, cf.~Section~\ref{sec:singheat} below),
but we just get a ``global'' (i.e., integrated)
version of the entropy inequality.

Moreover, as a further drawback of this procedure,
we note the explicit appearance of
the pressure $p$ in \eqref{i24},  which
explicitly
enters the system as an unknown.
To get a control of it from the momentum equation,
specific choices of the boundary conditions are needed
(cf.~\cite{MAL2} for further comments on this
point). In particular, this can be
reached in our case of periodic conditions.
These being said, we can state the

\greybox{
\centerline{\textsc{Total energy balance:}}

\vspace{-6mm}
\bFormula{i24}
  \partial_t \left( \frac{1}{2} |\vu|^2 + e \right)
  + \Div \left( \left( \frac{1}{2} |\vu|^2 + e \right) \vu \right)
  + \Div \vc{q}
  = \Div ( \sigma \vu )
   + \Div \big( \Gamma(\vt) \Grad \QQ : \HH \big),
\eF
}
\noindent
where
\begin{equation}\label{kappaciapo}
  \Div \vc{q}
   = - \Delta \kappaciapo(\vt)
    - \Div ( \Grad \QQ : \SSS (\Grad\vu, \QQ) ),
    \ \ \kappaciapo(\vt) = A_0 \vt + \frac{A_k}{k+1} \vt^{k+1}.
\end{equation}
In \eqref{i24}, we have noted as $e$
the {\sl internal energy}
\bFormula{defie}
  e = \calF + \vt s
   = \frac12 | \Grad \QQ |^2 - \frac\lambda2 |\QQ|^2 + \vt
    + a (m-1) \vt^m,
\eF
where the {\sl entropy} $s$ is given by
\bFormula{entropy}
  s = - \frac{\partial \mathcal{F}}{\partial \vt}
   = 1 + \log \vt - f(\QQ) + m a \vt^{m - 1}.
\eF
Recall also that the stress $\sigma$ is given by \eqref{i25}.

If one starts from the strong formulation, then
\eqref{i24} is (formally) obtained by multiplying
\eqref{i22} by $\vu$, \eqref{i23} by $\HH$, and
summing the results to \eqref{p6}. This procedure
involves a number of lengthy, but otherwise straighforward,
computations, which are provided in
detail in \cite[Section~1]{FRSZ}. In particular,
one uses the identity
\bFormula{mat}
 - \HH : \SSS(\Grad \vu, \QQ)
\eF
\[
 = \left( \QQ \HH - \HH \QQ \right) : \Grad \vu
  + 2 \xi \left( \HH: \QQ \right) \left( \QQ : \Grad \vu \right)
  - \xi \left[ \HH \left( \QQ  + \frac{1}{3} \II \right)
   + \left( \QQ  + \frac{1}{3} \II \right) \HH \right]: \Grad \vu,
\]
which holds for any symmetric matrix $\HH$.

Finally, multiplying the heat equation \eqref{p6}
 by $1/\vt$, one readily deduces the ``pointwise version''
of the entropy production inequality:
\bFormula{i28}
  \partial_t s
   + \Div ( s \vu )
   - \Div \left( \frac{\kappa(\vt)}{\vt} \Grad \vt \right)
\eF
\[
  \geq \frac{1}{\vt} \left( \frac{\mu(\vt)}{2} \left| \Grad \vu + \Grad^t \vu \right|^2
   + \Gamma(\vt) | \HH |^2
   + \frac{\kappa(\vt)}{\vt} |\Grad \vt |^2 \right).
\]
As noted above, in our existence result we are not able to
deal with relation \eqref{i28}, even interpreted as a
distributional inequality. Indeed, an $L^p$-estimate for
the term $s\vu$ on the \lhs\ seems out of reach due to
the occurrence of the singular function $f(\QQ)$ as a
summand in $s$ (cf.~\eqref{entropy}). Hence,  we have
to replace \eqref{i28}, by the following ``global''
entropy balance (simply obtained integrating
\eqref{i28} both in space and in time):

\medskip

\greybox{
\centerline{\textsc{Integrated entropy production inequality:}}

\vspace{-6mm}
\begin{equation}\label{i28weak}
   \ito\frac{1}{\vt} \left( \frac{\mu(\vt)}{2} \left| \Grad \vu + \Grad^t \vu \right|^2
    + \Gamma(\vt) | \HH |^2
    + \frac{\kappa(\vt)}{\vt} |\Grad \vt |^2 \right)
    - \io (s(t)-s(0))
   \le 0, \quad  \hbox{for a.e. }t\in (0,T).
\end{equation}
}
\noindent%

We are now ready to detail the notion of {\sl weak solution}\/
to our problem. Recall that the exponent $A$ was
introduced in \eqref{defiA}.
\bede\label{def:weak}
 A {\rm weak solution} to the non-isothermal liquid crystal
 system with the Ball-Majumdar potential is a {\rm quadruplet}
 $(\vu,p,\QQ,\vt)$ enjoying the regularity properties
 \begin{align}\label{rego:u}
   & \vu \in L^\infty(0,T;L^2(\Omega;\RR^3)),\ \
    \Grad\vu \in L^{\frac{2A}{A+1}}((0,T)\times\Omega;\RR^{3\times3}),\\
  \label{rego:p}
   & p \in L^{3/2}((0,T)\times\Omega),\\
  \label{rego:Q}
   & \QQ \in L^\infty(0,T;H^1(\Omega;\TS))\cap
    L^2(0,T;H^2(\Omega;\TS)), \ \
    \QQ_t \in L^{\frac{12}{7}}((0,T)\times\Omega;\TS),\\
  \label{rego:teta}
   & \vt \in L^A((0,T)\times \Omega), \ \
    \teta>0 \ \text{a.e.~in~$(0,T)\times \Omega$,}\\
  \label{rego:teta-1}
   & \log \vt \in L^\infty(0,T;L^1(\Omega)) \cap
    L^2(0,T;H^1(\Omega)),\\
\label{rego:toten}
&\frac12 |\vu|^2 + e \in C_w([0,T]; X), 
 \end{align}
where $X$ is a Sobolev space of negative order depending on $k$ and $m$, 
 satisfying the incompressibility condition~\eqref{i21},
 the momentum equation~\eqref{i22}, the Q-tensor
 evolution system~\eqref{i23} and the total energy
 balance~\eqref{i24}
 at least in the sense of distributions,
 and complying with the initial conditions in
 the following sense:
 \begin{align}\label{init:w0}
   & \vu|_{t=0} = \vu_0, \ \
   \QQ|_{t=0} = \QQ_0, \\
  \label{init:w}
   & \Big( \frac12 |\vu|^2 + e \Big)\bigg|_{t=0}
     = \Big( \frac12 |\vu_0|^2 + e_0 \Big),
 \end{align}
 where
 \bFormula{defie0}
   e_0 := \frac12 | \Grad \QQ_0 |^2
    - \frac\lambda2 |\QQ_0|^2
    + \vt_0 + a (m-1) \vt_0^m
 \eF
 and $e$ is defined in \eqref{defie}.
\edde
\noindent%
Then, our existence theorem, which can be considered
as the main result of this paper, reads

\greybox{
\bete\label{th:main}
 Under the assumptions stated in\/ {\rm Subsection~\ref{subsec:hyp}},
 the non-isothermal liquid crystal model with Ball-Majumdar
 potential admits at least one weak solution, which additionally
 satisfies the integrated entropy production
 inequality~\eqref{i28weak}.
\ente
}
\noindent%
The main steps of the proof of the theorem are carried out
in the next section.
\beos\label{weakprecise}
 In the statement, we claimed that the system equations
 hold at least in the sense of distributions. To
 be more precise, we can say that
 the momentum equation holds in the form
 \bFormula{m3}
   \int_0^T \io \Big[ \vu \cdot \partial_t \varphi + (\vu \otimes \vu) : \Grad \varphi \Big]
    = \int_0^T \io \sigma : \Grad \varphi
     - \io \vu_0 \cdot \varphi(0, \cdot),
 \eF
 satisfied for any test function $\varphi \in \DC([0,T) \times \Omega; R^3)$.
 %
 %
 %
 On the other hand, by virtue of the regularity properties
 \eqref{rego:u}-\eqref{rego:teta-1}, the order parameter
 equation \eqref{i23} turns out to
 hold pointwise (i.e.~in some $L^p$-space) and
 the periodic b.c.~for $\QQ$ hold in the sense of traces.
 Finally, the total energy balance (\ref{i24})
 is satisfied as the integral identity
 \bFormula{m5}
  \int_0^T \io \left[ \left( \frac{1}{2}|\vu|^2 + e \right) \partial_t \varphi
    + \left( \frac{1}{2}|\vu|^2 + e \right)
   \vu \cdot \Grad \varphi
   + \kappaciapo(\vt) \Delta \varphi
   - ( \Grad \QQ : \SSS(\Grad \vu, \QQ) ) \cdot \Grad \varphi \right]
 \eF
 \[
  = \int_0^T \io \Big[ \sigma \vu \cdot \Grad \varphi
     + \Gamma(\vt) (\Grad  \QQ : \HH) \cdot \Grad \varphi \Big]
  - \io \left( \frac{1}{2}|\vu_0|^2 + e_0 \right) \varphi(0, \cdot)
 \]
 for any $\varphi \in \DC([0,T) \times \Omega)$, where
 $e$ was defined in \eqref{defie} and $\kappaciapo$
 was introduced in \eqref{kappaciapo}. Actually, we needed
 to substitute the expression \eqref{i26} for $\vc{q}$
 and manage it by a further integration by parts.
\eddos

 Finally, let us note here that in Section~\ref{sec:singheat} we will state another existence result (Thm.~\ref{th:main:2}) holding in case of singular heat flux law \eqref{h12s} and in that case the entropy inequality holds true in the usual distributional form (cf. \eqref{i28we}).


\section{Proof of Theorem~\ref{th:main}}
\label{sec:proof}

The remainder of the paper is devoted to the proof of
Theorem~\ref{th:main}. The argument we are going to perform
follows closely the lines of the proof given in
our previous work~\cite{FRSZ}, the main difference being
represented of course by the presence of
the Ball-Majumdar potential.
In order to reduce complications, we will avoid
describing in detail the
(rather long and technical) approximation procedure
needed to get an existence theorem for a regularized
version of the system. Indeed, at this level, we do not see
relevant differences with respect to the argument
given in \cite{FRSZ}, to which we refer the reader
for more details.

Instead, the subsequent part of the proof
is more significantly different from that in \cite{FRSZ}
and, for this reason, will be presented in an extended
way. As a first step, we will prove in
Subsection~\ref{sec:apriori} that any
weak solution satisfies (at least formally)
a number of a-priori estimates, mainly coming
as direct consequences of the energy and entropy
balances. In particular, we will get a control
of certain norms of the solution components only
in terms of the initial data and of the fixed
parameters of the system.
These basic estimates, however, do not provide any
bound on the (subdifferential of the)
singular term $f$. Its control is actually
based on a rather delicate argument  detailed  in
Section~\ref{sec:BM} and concerning estimates of its Hessian.
This part is the key point of our argument and
represents the main novelty of the present paper.

In addition to that, we have to stress that,
compared to \cite{FRSZ}, we have here different
(power-like) growth conditions on the
viscosity coefficient $\Gamma(\vt)$, on the
heat conductivity (cf.~\eqref{h12}),
and on the specific heat (cf.~\eqref{i1}).
Mathematically speaking, these are needed in order
to control some $L^p$-norm of $\vt$.  However, in this case, we can allow 
for a general dependence of the viscosity coefficient from the 
absolute temperature (cf. \eqref{h11}). 

With the a-priori estimates at our disposal,
we will prove in Subsection~\ref{sec:limit} that
any sequence $(\vu_n,p_n,\QQ_n,\vt_n)$ which complies
with the bounds uniformly with respect
to~$n$ admits a limit point which is still a weak
solution of the system. This part, which we call
``weak sequential stability'' of weak solutions,
can be viewed as a simplified version of the compactness
argument needed to remove some regularization
of approximation of the system (like that provided
in \cite{FRSZ}).


\subsection{A priori estimates}
\label{sec:apriori}

In this section, we perform formal  {\it a priori} estimates.
As noted above, these estimates can be
rigorously justified for certain approximate
solutions in the way detailed in \cite{FRSZ}.
Most of the bounds follow directly
from the total energy conservation \eqref{i24}
and the entropy inequality \eqref{i28}
(indeed, the integrated version \eqref{i28weak},
holding for weak solutions, suffices for our
purposes).

\medskip

\noindent%
{\bf Energy-entropy estimate.}~~%
Integrating \eqref{i24} in space and using \eqref{defie} with the
periodic boundary conditions, we readily obtain
the {\sl energy estimate}
\begin{equation}\label{energy}
  \ddt \io \Big( \frac12 | \vu |^2 + e \Big)
   = \ddt \io \Big( \frac12 | \vu |^2 + \frac12 | \Grad \QQ |^2
     - \frac\lambda2 |\QQ|^2 + \vt + a (m-1) \vt^m \Big)
   = 0.
\end{equation}
Due to lack of coercivity of the internal
energy w.r.t.~$\QQ$,
the above information has to be coupled
with the bound coming from the
entropy inequality \eqref{i28weak}. Actually,
recalling \eqref{entropy} and \eqref{i28}, we infer
\bFormula{i28b}
  \ddt \io \big( - \log \vt - m a \vt^{m-1} + f(\QQ) \big)
   + \io \Big( \frac{\mu(\vt)}{2 \vt} \left| \Grad \vu + \Grad^t \vu \right|^2
   + \frac{\Gamma(\vt)}{\vt} | \HH |^2
   + \frac{\kappa(\vt)}{\vt^2} |\Grad \vt |^2 \Big) \le 0.
\eF
Summing \eqref{energy} to (some positive constant times)
\eqref{i28b}, we immediately obtain the estimates
\begin{align}\label{st11}
  & \| \vu \|_{L^\infty(0,T;L^2(\Omega;\RR^3))} \le c,\\
 \label{st12a}
  & \| \log \vt \|_{L^\infty(0,T;L^1(\Omega))} \le c,\\
 \label{st12}
  & \| \vt \|_{L^\infty(0,T;L^m(\Omega))} \le c,\\
 \label{st13}
  & \| f(\QQ) \|_{L^\infty(0,T;L^1(\Omega))} \le c,\\
 \label{st14}
  & \| \QQ \|_{L^\infty(0,T;H^1(\Omega))} \le c.
\end{align}
In particular, we used here the singular character of
$f$, which also gives
\begin{equation}\label{st18}
  \| \QQ \|_{L^\infty((0,T)\times \Omega;\TS)} \le c.
\end{equation}
Moreover, we have the information coming from the
latter integral in \eqref{i28b}, which we now specify.

First of all, using \eqref{h12}, we readily infer
\begin{equation}\label{st15}
  \| \Grad \vt^{\kappa/2} \|_{L^2(0,T;L^2(\Omega;R^3))}
   + \| \Grad \log\vt \|_{L^2(0,T;L^2(\Omega;R^3))}
   \le c.
\end{equation}
%
%

Combining \eqref{st12a}, the second part of \eqref{st15}
and using some generalized version of Poincar\'e's inequality
(see, e.g., \cite[Lemma~3.2]{GSZ}), we get
\begin{equation}\label{st15b}
  \| \log\vt \|_{L^2(0,T;H^1(\Omega))} \le c.
\end{equation}
Next, by \eqref{st12}, \eqref{st15} and interpolation,
recalling also \eqref{defiA}, it is easy to infer
\begin{equation}\label{st19}
  \| \vt \|_{L^A((0,T)\times \Omega)}
   = \| \vt \|_{L^{\frac{3k+2m}3}((0,T)\times \Omega)} \le c.
\end{equation}
Next, by virtue of \eqref{h11} and \eqref{h13}, we obtain
\begin{align}\label{st16}
  & \| \vt^{-1/2} \Grad \vu \|_{L^2(0,T;L^2(\Omega;\RR^{3\times 3}))} \le c,\\
 \label{st17}
  & \| \HH \|_{L^2(0,T;L^2(\Omega;\TS))}
   + \| \vt^{-1/2} \HH \|_{L^2(0,T;L^2(\Omega;\TS))} \le c.
\end{align}
As a consequence, by \eqref{st19}, \eqref{st16} and interpolation,
we have
\begin{equation}\label{st20}
  \| \Grad \vu \|_{ L^{\frac{2A}{A+1}}((0,T)\times \Omega;\RR^{3\times 3}) }  \le c,
\end{equation}
whence, also,
\begin{equation}\label{st20b}
  \| \vu \|_{ L^{\frac{2A}{A+1}}(0,T;L^{\frac{6A}{A+3}}(\Omega;\RR^3))}
   \le c \| \vu \|_{ L^{\frac{2A}{A+1}}(0,T;W^{1,\frac{2A}{A+1}}(\Omega;\RR^3))}
   \le c.
\end{equation}
Interpolating between \eqref{st11} and \eqref{st20b}
and using that $A>9$ (cf.~\eqref{defiA}), we arrive at
\begin{equation}\label{li11}
  \big\| | \vu |^3 \big\|_{L^{1+}((0,T)\times \Omega)} \le c.
\end{equation}
Here $1+$ denotes a generic exponent strictly greater than $1$. This
kind of notation will be used also in the sequel.
\beos\label{on:positivity}
 It is worth stressing that in all the procedure
 (cf.~for instance estimates~\eqref{st12} and \eqref{st15})
 we are implicitly assuming that $\vt$ is positive
 at least almost everywhere. Actually, in \cite{FRSZ}
 we were able to prove this property (and, actually,
 also a stronger one) by means of maximum principle arguments,
 even at the approximated level. In the current setting,
 due to the presence of singular terms on the \rhs\ of
 \eqref{p6}, ensuring the same property for \emph{approximate} solutions is more delicate, requiring additional 
 singular terms to be introduced in the process of approximation, see e.g. \cite[Chapter 3]{FEINOV}.
\eddos

\smallskip

\noindent%
{\bf Estimate of the singular potential.}~~%
%
%

The estimate of the (subdifferential
of the) Ball-Majumdar potential 
constitutes the key point of our argument. Recalling estimate \eqref{st17}, we compute the
$L^2$-scalar product of $\HH$ with
$- \vt \calL \left[ \frac{\partial f(\QQ)}{\partial \QQ } \right]$.
The necessary key properties of $f$  will be proved in the subsequent
Section~\ref{sec:BM}. Hence, we have
\begin{align}\label{test11}
  & I_1 + I_2
   := \io \vt^2 \calL \left[ \frac{\partial f(\QQ)}{\partial \QQ } \right]
      : \calL \left[ \frac{\partial f (\QQ)}{\partial \QQ } \right]
   - \io \vt \Delta \QQ : \calL \left[ \frac{\partial f (\QQ)}{\partial \QQ } \right]\\
  \nonumber
   & \mbox{}~~~~~
    = - \io \vt \HH : \calL \left[ \frac{\partial f(\QQ)}{\partial \QQ } \right]
     + \lambda \io \vt \QQ : \calL \left[ \frac{\partial f(\QQ)}{\partial \QQ } \right]
     =: J_1 + J_2.
\end{align}
%

%
%
Moreover,
by H\"older's and Young's inequalities and uniform boundedness
of $\QQ$ (cf.~\eqref{st18}),
\begin{equation}\label{test13}
  J_2 = \lambda \io \vt \QQ : \calL \left[ \frac{\partial f(\QQ)}{\partial \QQ } \right]
   \le \frac\kappa4 \left\| \vt \calL \left[ \frac{\partial f(\QQ)}{\partial \QQ } \right]
      \right\|_{L^2(\Omega;\TS)}^2 + c(\kappa), \quad \kappa > 0
\end{equation}
and
\begin{equation}\label{test14}
  J_1 = - \io \vt \HH : \calL \left[ \frac{\partial f(\QQ)}{\partial \QQ } \right]
   \le \frac\kappa4 \left\| \vt \calL \left[ \frac{\partial f(\QQ)}{\partial \QQ } \right]
      \right\|_{L^2(\Omega;\TS)}^2
    + c(\kappa) \| \HH \|_{L^2(\Omega;\TS)}^2.
\end{equation}
Notice that, integrated in time, the last term can be estimated
by virtue of the first one in \eqref{st17}.

The key point is to control the term $I_2$ on the \lhs\ of
\eqref{test11}. Using Green's formula, we have
\begin{align}\label{test15}
  I_2 & = - \io \vt \Delta \QQ : \calL \left[ \frac{\partial f(\QQ)}{\partial \QQ } \right]\\
  \nonumber
   & = \io \vt \frac{\partial^2 f(\QQ)}{\partial \QQ_{ij} \partial \QQ_{kl} } \Grad \QQ_{ij}
   \cdot \Grad \QQ_{kl}
   + \io \left( \calL \left[ \frac{\partial f(\QQ)}{\partial \QQ } \right] : \Grad \QQ \right) \cdot \Grad \vt
    =: I_{2,1} + I_{2,2}.
\end{align}
Then, thanks to the key estimate \eqref{f:ftest1} below and  \eqref{st18},
\begin{equation}\label{test15b}
  I_{2,1}
   = \io \vt \frac{\partial^2 f(\QQ)}{\partial \QQ_{ij} \partial \QQ_{kl} } \Grad \QQ_{ij}
               \cdot \Grad \QQ_{kl}
   \ge \ep \io \vt \left| \calL \left[ \frac{\partial f(\QQ)}{\partial \QQ } \right] : \Grad \QQ \right|^2.
\end{equation}
On the other hand, by Young's inequality,
\begin{equation}\label{test16}
  | I_{2,2} | \le \frac\ep2
       \io \vt \left| \calL \left[ \frac{\partial f(\QQ)}{\partial \QQ } \right]
            : \Grad \QQ \right|^2
   + c(\ep) \io \left| \frac{\Grad\vt}{\vt^{1/2}} \right|^2.
\end{equation}
Then, integrating \eqref{test11} over $(0,T)$,  and taking $\kappa=\ep/2$,   we get 
%
\begin{equation}\label{st21a}
  \left\| \vt \calL \left[ \frac{\partial f(\QQ)}{\partial \QQ } \right] \right\|_{L^2((0,T)\times \Omega;\TS)}
    \le c,
\end{equation}
whence, comparing terms in $\HH$ and using \eqref{st18} and
the first \eqref{st17}, we arrive at
\begin{equation}\label{st21}
  \| \Delta \QQ \|_{L^2((0,T)\times \Omega;\TS)}
    \le c.
\end{equation}

Summarizing, we have obtained an $L^2$-control on both the Laplacian of $\QQ$ and the
singular term in \eqref{i23}.


\subsection{Weak sequential stability}
\label{sec:limit}

In this part, we assume to have a sequence $(\vu_n,p_n,\QQ_n,\vt_n)$
of weak solutions satisfying the estimates proved in the
previous subsection uniformly with respect to $n$.
Under these conditions, we show that there exists 
subsequence that converges in a suitable way to another
weak solution to the problem. In order to simplify
the notation, we assume all convergence relations
appearing in the sequel to hold up to the
extraction of (non-relabelled) subsequences.

First of all, notice that, thanks to
the bounds \eqref{st11}-\eqref{li11},
we have
\begin{align}\label{lim11}
  & \vu_n \to \vu \ \ \text{weakly star in }
   L^\infty(0,T;L^2(\Omega;\RR^3))
   \cap L^{\frac{2A}{A+1}}(0,T;W^{1,\frac{2A}{A+1}}(\Omega;\RR^3)), \\
 \label{lim12}
  & \QQ_n \to \QQ \ \ \text{weakly in }
   L^2(0,T;H^2(\Omega;\TS)),\\
 \label{lim12b}
  & \QQ_n \to \QQ \ \ \text{weakly star in }
   L^\infty(0,T;H^1(\Omega;\TS))
   \cap L^\infty(0,T;L^\infty(\Omega;\TS)), \\
 \label{lim13}
  & \vt_n \to \vt \ \ \text{weakly star in }
   L^\infty(0,T;L^m(\Omega))
   \cap L^A((0,T)\times \Omega)
   \cap L^2(0,T;H^1(\Omega)),
\end{align}
the latter relation coming from \eqref{st15} and
interpolation.
As a consequence of \eqref{st17}, we also obtain
\begin{equation}\label{lim14}
  \HH_n \to \ov{\HH} \ \
   \text{weakly in } L^2(0,T;L^2(\Omega;\TS)),
\end{equation}
for some tensor-valued function $\ov{\HH}$.
Here and below we use the upper bar to indicate weak limits
that are not yet identified. In particular,
we cannot express $\ov{\HH}$ in terms of
$\vt$ and $\QQ$ at this level.

Next, by \eqref{lim12}-\eqref{lim12b}
and Gagliardo-Nirenberg inequalities
(cf., e.g., \cite{nier}),
\begin{equation}\label{lim16}
  \QQ_n \to \QQ \ \
   \text{weakly in } L^4(0,T;W^{1,4}(\Omega;\TS)).
\end{equation}
Then, by \eqref{i5}, \eqref{st20}, $A>9$ (cf. \eqref{defiA}), and uniform
boundedness of $\QQ_n$,
\begin{equation}\label{lim17}
  \| \SSS(\Grad\vu_n,\QQ_n) \|_{L^{\frac{9}{5}+}((0,T)\times\Omega;\TS)}
   \le c,
\end{equation}
whereas, by \eqref{li11} and \eqref{lim16},
\begin{equation}\label{lim18}
  \| \vu_n \cdot \Grad \QQ_n \|_{L^{\frac{12}{7}+}((0,T)\times\Omega;\TS)}
   \le c.
\end{equation}
By \eqref{lim15}, \eqref{lim17}, \eqref{lim18}
and a comparison of terms in \eqref{i23}, we then infer
\begin{equation}\label{lim19}
  \partial_t \QQ_n \to \partial_t \QQ
   \ \ \text{weakly in, say, } L^{\min\{\frac{12}{7}+, 2A/(A+2)\}}((0,T)\times\Omega;\TS).
\end{equation}
Hence, by the Aubin-Lions lemma we get the
following strong convergence relations:
\begin{align}\label{lim19b}
  & \QQ_n \to \QQ \ \ \text{in $L^p((0,T)\times\Omega;\TS)$, for all $p\in [1,+\infty)$},\\
 \label{lim19b2}
  & \Grad\QQ_n \to \Grad\QQ \ \ \text{in $L^p((0,T)\times\Omega;R^{27})$,
    for all $p\in [1,4)$}.
\end{align}
%
%

\smallskip

To proceed, we now compute the limit $n\nearrow\infty$ in the
momentum equation. To start with, we first deal with the
pressure~$p$. This enters the total energy balance \eqref{i24}
as an additional unknown; hence, it needs to be estimated
directly. To this purpose, one (formally) applies the operator
$\Div$ to the momentum equation \eqref{i22}. This gives
rise to an elliptic problem for $p$ which is well-posed
thanks to the choice of periodic boundary conditions
for all unknowns  (cf. \cite{BFM} for more details on this point). Recalling \eqref{i25} and noting that
\begin{align}\label{lim20a}
  & \| \vu_n \otimes \vu_n \|_{L^{\frac32+}((0,T)\times\Omega;R^{3\times3})}
   + \| \mu(\vt_n) \Grad \vu_n \|_{L^{\frac95+}((0,T)\times\Omega;R^{3\times3})} \le c, \\
 \label{lim20b}
  & \| \HH_n \|_{L^2((0,T)\times\Omega;\TS)}
   + \| \Grad \QQ_n \odot \Grad \QQ_n \|_{L^2((0,T)\times\Omega;R^{3\times3})} \le c,
\end{align}
thanks to \eqref{li11}, \eqref{lim11}, \eqref{lim14}
and \eqref{lim16}, using also the uniform boundedness
of $\QQ_n$.  It is then not difficult to deduce
(cf.~\cite[Section~4]{FFRS} for more details)
\begin{equation}\label{lim21}
  p_n \to p \ \
   \text{weakly in } L^{\frac32+}((0,T)\times\Omega).
\end{equation}
Hence, a comparison of terms in \eqref{i22} gives that
\begin{equation}\label{lim22}
  \| \partial_t \vu_n \|_{L^{\frac32+}(0,T;X)}
   \le c,
\end{equation}
where $X$ is a suitable Sobolev space of negative order.
By the Aubin-Lions lemma and \eqref{li11} we then get
\begin{equation}\label{lim23}
  \vu_n \to \vu \ \
   \text{strongly in } L^{3+}((0,T)\times\Omega;R^3).
\end{equation}

\smallskip

Hence, we have deduced strong (and consequently pointwise)
convergence of $\QQ_n$ and $\vu_n$. The next step
consists in obtaining an analogous property
for the temperature. To this
aim, we consider the total energy balance at the level $n$ and
notice that
\begin{align}\label{lim31}
  & \| \Grad\QQ_n : \HH_n \|_{L^{\frac43}((0,T)\times\Omega;R^3)} \le c
   \ \ \text{by \eqref{lim14} and \eqref{lim16}},\\
 \label{lim32}
  & \| \sigma_n \vu_n \|_{L^{1+}((0,T)\times\Omega;R^3)} \le c
   \ \ \text{by \eqref{lim11}, \eqref{lim21} and \eqref{lim23}},\\
 \label{lim33}
  & \| |\vu_n|^2 \vu_n \|_{L^{1+}((0,T)\times\Omega);R^3)} \le c
   \ \ \text{by \eqref{lim23}}.
\end{align}
%
%
%
Using \eqref{lim16}, \eqref{i5}, \eqref{lim11} and uniform
boundedness of $\QQ_n$, we also get
\begin{equation}\label{lim35}
  \big\| \Grad \QQ_n : \SSS (\Grad\vu_n, \QQ_n) \big\|%
    _{ L^{\frac{36}{29}}((0,T)\times\Omega;R^3) } \le c.
\end{equation}
Moreover, using \eqref{lim13} and \eqref{defiA},
we have
\begin{equation}\label{lim36}
  \| \kappaciapo(\vt_n) \|_{ L^{1+}((0,T)\times\Omega)} \le c,
   \ \ \text{provided that} \ A=\frac{3k+2m}3 > k+1, \
   \text{i.e.,} \ m>\frac32,
\end{equation}
the latter condition being a consequence of \eqref{defiA}.
Recall that $\kappaciapo$ was defined in \eqref{kappaciapo}.

Now, it is clear that all summands in the internal energy~$e_n$
(cf.~\eqref{defie}) are uniformly $L^2$-bounded, with the
exception of the last one for which we have
\begin{equation}\label{lim37}
  \| \vt_n^m \|_{ L^{\frac32}((0,T)\times\Omega)} \le c
  \ \ \text{provided that} \ \frac{3k+2m}{3} \ge \frac{3m}2,
  \ \text{i.e.,} \ 6k \ge 5m,
\end{equation}
which also follows from \eqref{defiA}.
By \eqref{lim37} and \eqref{li11}, we then find
\begin{equation}\label{lim38}
  \| e_n \vu_n \|_{ L^{1+}((0,T)\times\Omega;R^3)} \le c.
\end{equation}
%
%
Collecting \eqref{lim31}-\eqref{lim38} and comparing terms
in \eqref{i24}, we then infer
\begin{equation}\label{lim39}
  \left\| \partial_t \left( \frac12 |\vu_n|^2 + e_n \right) \right\|_{L^{1+}(0,T;X)}
   \le c,
\end{equation}
where $X$ is, again, some Sobolev space of negative order.
%
%
On the other hand, a direct computation based on
\eqref{st15}, \eqref{st19} and the other estimates
permits to check that
\begin{equation}\label{lim41}
  \left\| \Grad \left( \frac12 |\vu_n|^2 + e_n \right) \right\|_{ L^{1+}((0,T)\times\Omega;R^3)} \le c,
  \ \ \text{provided that} \ \frac{3k+2m}{3} > 2m-k,
  \ \text{i.e.,} \ 3k > 2m,
\end{equation}
and also this condition  is satisfied due to \eqref{defiA}.
Actually, we used here that
\begin{equation}\label{lim42}
  \Grad \vt^m = c_{k,m} \vt^{m-\frac{k}2} \Grad \vt^{\frac{k}2}.
\end{equation}
With \eqref{lim39} and \eqref{lim41} at disposal, we can use
once more the Aubin-Lions lemma to conclude that
\begin{equation}\label{lim43}
  \left( \frac12 |\vu_n|^2 + e_n \right) \to \left( \frac12 |\vu|^2 + e \right) \ \
   \text{strongly in } L^{1+}((0,T)\times\Omega;R^3)
\end{equation}
and, consequently, almost everywhere in $(0,T)\times \Omega$.
Note in particular that the limit of $|\vu_n|^2$ is identified
as $|\vu|^2$ thanks to \eqref{lim23}. Moreover,
by \eqref{lim11}, \eqref{lim13}, and \eqref{lim12b},
we also have
\begin{equation}\label{lim43b}
  \left\| \frac12 |\vu_n|^2 + e_n \right\|_{L^\infty(0,T;L^1(\Omega))}
   \le c.
\end{equation}
Hence, \eqref{lim39}, \eqref{lim43b} and a generalized
form of the Aubin-Lions lemma give
\begin{equation}\label{lim43c}
  \left( \frac12 |\vu_n|^2 + e_n \right) \to \left( \frac12 |\vu|^2 + e \right) \ \
   \text{strongly in } C^0([0,T];X),
\end{equation}
where $X$ is, again, some Sobolev space of negative order. Hence, the
Cauchy condition for the total energy balance holds also
in the limit, in the sense specified by \eqref{init:w}.

By \eqref{lim43} and \eqref{lim23}, we also obtain strong
convergence of $e_n$, say, in $L^{1+}$. Consequently, computing
\begin{equation}\label{lim44}
  \int_0^T \io (e_n - e_m) \sign (\vt_n - \vt_m)
\end{equation}
for a couple of indexes $n$, $m$, and using monotonicity of
the function $\vt \mapsto \vt^m$ (recall that the temperature
is assumed to be nonnegative) and strong $L^2$-convergence of
$\QQ_n$ and $\Grad \QQ_n$, we readily obtain that $\{\vt_n\}$ is
a Cauchy sequence in $L^1((0,T)\times \Omega)$. Hence, recalling
\eqref{lim13} and using a proper generalized version of
Lebesgue's theorem, we arrive at
\begin{equation}\label{lim45}
  \vt_n \to \vt \ \
   \text{strongly in} \ L^p((0,T)\times \Omega)
   \ \text{for all }p\in [1,A).
\end{equation}
As a consequence, thanks to assumptions \eqref{h11}-\eqref{h13}, we infer
\begin{align}\label{limh11}
  & \mu(\vt_n) \to \mu(\vt) \ \
   \text{strongly in} \ L^p((0,T)\times \Omega)
   \ \text{for all }p\in [1,\infty), \\
 \label{limh12}
  & \kappa(\vt_n) \to \kappa(\vt) \ \
   \text{strongly in} \ L^p((0,T)\times \Omega)
   \ \text{for all }p\in \Big[1,\frac{3k+2m}{3k}\Big),\\
 \label{limh13}
  & \Gamma(\vt_n) \to \Gamma(\vt) \ \
   \text{strongly in} \ L^p((0,T)\times \Omega)
   \ \text{for all }p\in [1,A).
\end{align}
Combining \eqref{lim13} and \eqref{lim14}, and recalling
assumption \eqref{h13}, we also get
 \begin{equation}\label{lim15}
  \Gamma(\vt_n)\HH_n \to \ov{\Gamma(\vt)\ov{\HH}} \ \
   \text{weakly in } L^{\frac{2A}{A+2}}(0,T;L^{\frac{2A}{A+2}}(\Omega;\TS)).
\end{equation}
Moreover, relation \eqref{limh13} also implies
$\ov{\Gamma(\vt)\ov{\HH}} = \Gamma(\vt)\ov{\HH}$  (cf. also \eqref{lim15}).
This permits us to pass to the limit in the momentum
equation \eqref{i22} (in particular, the stress
$\sigma$ is identified in terms of the limit
functions $\vt$, $\QQ$, $\ov{\HH}$ and $\vu$). Analogously,
we can take the limit in the Q-tensor equation
\eqref{i23} (where, however, the function $\ov{\HH}$ is still
to be identified) and in the total energy balance
\eqref{i24} (we use here also properties
\eqref{lim35}-\eqref{lim39}).

In addition, from \eqref{lim45} and \eqref{st15b} we also
obtain
\begin{equation}\label{lim45b}
  \log \vt_n \to \log \vt \ \
   \text{strongly in} \ L^p((0,T)\times \Omega)
   \ \text{for all }p\in [1,2).
\end{equation}
In particular, if positivity holds (almost everywhere) for
$\vt_n$, then it is preserved in the limit $\vt$.

To conclude the proof, we need to take $n\nearrow +\infty$
in the relations involving the singular potential $f(\QQ)$.
To be precise, what remains to do is identifying
$\ov{\HH}$ in terms of $\vt$ and $\QQ$ and letting
$n\nearrow\infty$ in the weak form \eqref{i28weak}
of the entropy production inequality.

We start with the first task.
Thanks to estimate \eqref{st21a}, we get that
\begin{equation}\label{lim51}
  \vt_n \calL \left[ \frac{\partial f(\QQ_n)}{\partial \QQ_n } \right]
   \to \ov{\vt \calL \left[ \frac{\partial f(\QQ)}{\partial \QQ } \right]} \ \
  \text{weakly in} \ L^2((0,T)\times \Omega;\TS).
\end{equation}
Moreover, by strong convergence of
$\vt_n$ and $\QQ_n$ (cf.~\eqref{lim45} and \eqref{lim19b}),
and using that both $\QQ_n$ and $\QQ$ take
their values almost everywhere into $\calD[f]$ (otherwise $f(\QQ_n)$ could
not be uniformly $L^1$-bounded), we obtain
\begin{equation}\label{lim52}
  \vt_n \calL \left[ \frac{\partial f(\QQ_n)}{\partial \QQ_n } \right]
   \to \vt \calL \left[ \frac{\partial f(\QQ)}{\partial \QQ } \right] \ \
  \text{a.e.~in} \ (0,T)\times \Omega.
\end{equation}
By a generalized form of Lebesgue's Theorem, the combination
of these facts entails that
\begin{equation}\label{lim53}
  \vt_n \calL \left[ \frac{\partial f(\QQ_n)}{\partial \QQ_n } \right]
   \to \vt \calL \left[ \frac{\partial f(\QQ)}{\partial \QQ } \right] \ \
  \text{strongly in} \ L^p((0,T)\times \Omega;\TS) \ \
  \text{for all}\ p\in [1,2).
\end{equation}
As a consequence of \eqref{lim53}, $\ov{\HH}$ is identified
as $\HH$, in terms of $\vt$ and $\QQ$; namely, \eqref{i20} holds.
In particular, for $n\nearrow \infty$,
equation \eqref{i23} goes to the expected limit.

To conclude the proof, we have to take $n\nearrow\infty$ in
the ``weak'' entropy production inequality \eqref{i28weak}.
To this aim, we notice that, for any nonnegative-valued
$\phi \in C_c^\infty([0,T)\times \Omega)$, there holds
\begin{align}\label{i28lim}
  & \iTo \frac{\phi}{\vt} \left( \frac{\mu(\vt)}{2} \left| \Grad \vu + \Grad^t \vu \right|^2
   + \Gamma(\vt) | \HH |^2
   + \frac{\kappa(\vt)}{\vt} |\Grad \vt |^2 \right)\\
 \nonumber
  & \mbox{}~~~~~
   \le \liminf_{n\nearrow\infty}
     \iTo \frac{\phi}{\vt_n} \left( \frac{\mu(\vt_n)}{2} \left| \Grad \vu_n + \Grad^t \vu_n \right|^2
   + \Gamma(\vt_n) | \HH_n |^2
   + \frac{\kappa(\vt_n)}{\vt_n} |\Grad \vt_n |^2 \right).
\end{align}
This is, indeed, a consequence of relations \eqref{lim11},
\eqref{lim13}, \eqref{lim14} (where, now, $\ov{\HH}=\HH$),
of convexity of the above integrand with respect to
$\Grad \vu + \Grad^t \vu$, $\HH$, and $\Grad \vt$,
and of a standard semicontinuity argument
(see, e.g., \cite{Ioffe}). Choosing $\phi\equiv1$,
we can manage the first integral in \eqref{i28weak}.
%
%
Next, in order to deal with the second integral of
\eqref{i28weak}, we simply observe that
$-s$ is  the sum of convex function both of $\vt$ and of $\QQ$ and of the function 
$ma\vt^{m-1}$
(cf.~\eqref{entropy}). Hence, we can take the
(infimum) limit of that integral for a.e. $t\in (0,T)$  by virtue of
relations (\ref{lim12})-(\ref{lim13}) and \eqref{lim45},
and of a further semicontinuity argument. This
concludes the proof of Theorem \ref{th:main}.



\section{Analytical properties of the Ball-Majumdar potential}
\label{sec:BM}

In this part we prove a number of fine properties of the Ball-Majumdar
potential that were used in the proof.  We recall first the definition and the basic features of the Ball-Majumdar potential. 
 As Observed in \cite{ball-lecture, BalMaj},the Ball-Majumdar singular potential $f$, defined as a function of $Q$, can be expressed  as:

\be
f(Q)=F^{BM}(\lambda_1,\lambda_2,\lambda_3)
\ee  where $(\lambda_1,\lambda_2,\lambda_3)$ are the eigenvalues of $Q$ and $F^{BM}$ is defined as

\begin{equation}\label{def:FBM}
F^{BM}(\lambda_1,\lambda_2,\lambda_3)=\sum_{i=1}^3\mu_i (\lambda_i+\frac{1}{3})-\ln Z(\mu_1,\mu_2,\mu_3)
\end{equation} with

\be\label{def:Z}
Z(\nu_1,\nu_2,\nu_3)=\int_{\Sphere^2} \exp(\sum_{j=1}^3\nu_j p_j^2)\, {\rm d} p.
\ee and $\mu_i,i=1,2,3$, are given implicitly as solutions of the system:

\begin{equation}\label{def:musystem}
\frac{\partial \ln Z}{\partial\mu_i}=\lambda_i+\frac{1}{3},\,i=1,2,3.
\end{equation}

The solutions of \eqref{def:musystem} are determined up to an additive constant added to all the $\mu_i$ (i.e. if $(\mu_1,\mu_2,\mu_3)$ is a solution for a given
triplet $(\lambda_1,\lambda_2,\lambda_3)$ then so is $(\mu_1+C,\mu_2+C,\mu_3+C)$ for any $C\in\RR$).

Nevertheless, let us note that the $\mu_i,i=1,2,3$, can be taken, without loss of generality, as real analytic functions of $(\lambda_1,\lambda_2,\lambda_3)$. More precisely we have:

\begin{lemm}\label{lemma:mudiff}
There exist $\mu_i,i=1,2,3$, real analytic functions of $(\lambda_1,\lambda_2,\lambda_3)$, such that
\begin{equation}\label{mu:sumconstraint}
\sum_{i=1}^3 \mu_i=0
\end{equation} and solving the system \eqref{def:musystem}.
\end{lemm}
\beos
The condition \eqref{mu:sumconstraint} is just a convenient way of eliminating the one-dimensional degeneracy in the solutions of the system \eqref{def:musystem}.
\eddos
\bproof
We consider the linear space
\be\label{def:X}
X=\{y\in \RR^3: y\cdot (1,1,1)=0\}.
\ee

 We aim to apply the analytic implicit function theorem to $\mathcal{F}:X\times X\to X$ defined by $$\mathcal{F}(\lambda_1,\lambda_2,\lambda_3,\mu_1,\mu_2,\mu_3)=\Big(\frac{\partial \ln Z}{\partial\mu_1}-\lambda_1-\frac{1}{3},\frac{\partial \ln Z}{\partial\mu_2}-\lambda_2-\frac{1}{3},\frac{\partial \ln Z}{\partial\mu_3}-\lambda_3-\frac{1}{3}\Big)$$ in order to express $(\mu_i)_{i=1,2,3}$ as functions of $(\lambda_i)_{i=1,2,3}$.

We note that $X$ is a $2$-dimensional space, with an orthonormal basis $f_1=\frac{1}{\sqrt{2}}(-1,1,0),f_2=\frac{1}{\sqrt{6}}(1,1,-2)$. As such, in terms of this orthonormal basis, the non-degeneracy condition needed for applying the implicit function theorem is $\det M\not=0$ for

$$M=\left\{\begin{array}{ll} \frac{\partial (\mathcal{F},f_1)}{\partial f_1} & \frac{\partial (\mathcal{F},f_1)}{\partial f_2}\\ \frac{\partial (\mathcal{F},f_2)}{\partial f_1} &\frac{\partial (\mathcal{F},f_2)}{\partial f_2}\end{array}\right\}.$$

One can check that $\det M=36\frac{\partial^2 \ln Z}{\partial \mu_1^2}\frac{\partial^2 \ln Z}{\partial\mu_2^2}-36\left(\frac{\partial^2 \ln Z}{\partial \mu_1\mu_2}\right)^2$. It is convenient to consider the matrix

$$N:=\left\{\begin{array}{ll} \frac{\partial^2 \ln Z}{\partial \mu_1^2} &\frac{\partial^2 \ln Z}{\partial\mu_1\partial\mu_2}\\ \frac{\partial^2 \ln Z}{\partial\mu_1\partial\mu_2} &\frac{\partial^2 \ln Z}{\partial \mu_1^2} \end{array}\right\}$$ and observe that $\det M=36\det N$. Then in order to check that $\det M\not=0$ it suffices to check that for $e=(e_1,e_2)\in \mathbb{R}^2$ we have $(Ne,e)\ge 0$  with equality if and only if $e=(0,0)$.
Indeed, taking into account the definition \eqref{def:Z} of $Z$ we have:

$$(Ne,e)=Z^{-2}\left[\left(\int_{\mathbb{S}^2} e^{\mu_k p_k^2}(p_\alpha^2 e_\alpha)^2\,dp\right)\left(\int_{\mathbb{S}^2} e^{\mu_k p_k^2}\,dp\right)-\left(\int_{\mathbb{S}^2} e^{\mu_k p_k^2}(p_\alpha^2e_\alpha)\,dp\right)^2\right]\ge 0$$ where for the last inequality we used Cauchy-Schwarz, and $k=1,2,3$, $\alpha=1,2$. It is clear that equality cannot hold for $e\not=(0,0)$ hence our claim.

\eproof

We can now prove the fundamental estimate on the Hessian of $f$:

\bepr\label{teo:ftest}
There  exists $\varepsilon>0$ such that:

 \begin{equation}\label{f:ftest1}
    \frac{\partial^2 f(\QQ)}{\partial \QQ_{ij} \partial \QQ_{kl} }
         \VV_{ij} \cdot \VV_{kl}
     \ge \varepsilon \left| \calL \left[ \frac{\partial f(\QQ)}{\partial \QQ } \right] :
       \VV \right|^2
     \ \
       \text{for all }\QQ \in \calD[f],\ \VV \in \TS.
 \end{equation}

\empr

\noindent  \bproof We proceed in several steps:

\bigskip\noindent{\it Step 1: reduction to proving the concavity of an auxiliary function}
\smallskip

We set $h(\QQ):=e^{-\eps f(\QQ)}$ for all $\QQ\in  \calD[f]$, where $\ep>0$ is a positive constant to be chosen later on. It is more convenient to work with $h$ because the required property \eqref{f:ftest1} is implied by $h$ being concave and the concavity of an isotropic function is in general simpler to check,  as it will be seen in the next step.

We claim thus that $h$ being concave, smooth on $\calD[f]$ and positive implies \eqref{f:ftest1}. We start by noting that for $\QQ \in\calD[f]$ and arbitrary $V\in\TS$:

\begin{align}\no
\varepsilon\frac{\partial f(\QQ)}{\partial \QQ_{ij}}:V_{ij}&=-\frac{\rm d}{\rm dt}\Big(\log h(\QQ+tV)\Big)|_{t=0}=-h(\QQ)^{-1}\left(\frac{\rm d}{\rm dt}h(\QQ+tV)\right)|_{t=0}\non\\
\no
&=- h(\QQ)^{-1}g(\QQ,V)
\end{align} where we denoted $g(\QQ,V):=\frac{\rm d}{\rm dt}h(\QQ+tV)|_{t=0}.$

Similarly, again for $\QQ\in\calD[f]$,

$$
\varepsilon\frac{\partial^2 f}{\partial \QQ_{ij}\partial \QQ_{mn}}V_{ij}V_{mn}:=-\frac{\rm d^2}{\rm dt^2}\log(h(\QQ+tV))|_{t=0}= \frac{-k(\QQ,V)h(\QQ)+g(\QQ,V)^2}{(h(\QQ))^2},
$$ where we denoted $k(\QQ,V):=\frac{\rm d^2}{\rm dt^2}h(\QQ+tV)|_{t=0}$.

Then the desired estimate \eqref{f:ftest1} is equivalent to
$$
\frac{g(\QQ,V)^2-k(\QQ,V)h(\QQ)}{(h(\QQ))^2}\ge\frac{g(\QQ,V)^2}{(h(\QQ))^2}
$$ which clearly holds as $h>0$  and $k:=\frac{\rm d^2}{\rm dt^2}h(\QQ+tV)|_{t=0}\le 0$ by our assumptions that $h$ is positive, smooth and concave.

The assumed positivity of $h$ is obvious (by the definition of $h$) and the smoothness is a consequence of the smoothness of $f$. We are left with checking the concavity.

\bigskip\noindent{\it Step 2: checking concavity for a function of eigenvalues}

\smallskip
We define now  $H:(-\frac{1}{3},\frac{2}{3})\to [0,\infty)$  by $H(\lambda_1,\lambda_2,\lambda_3):=h(\QQ)$ with $\QQ$ having eigenvalues $\lambda_i,\,i=1,2,3$ (it is known that such a function is well defined and symmetric, see for instance \cite{ball-isotropic} and the references therein). Then $h$ is concave if and only if $H$ is concave (see for instance Prop $18.2.4$ in \cite{silhavy}). We note that $h$ and $H$ have the same regularity (see for instance \cite{ball-isotropic}, Thm $5.5$ for up to $C^2$ regularity and the extension up to $C^\infty$ in \cite{silhavy-infty}). Thus $H$ is concave if and only if its Hessian is non-positive.

\bigskip\noindent{\it Step 3: the concavity of $H$ and the asymptotics of certain integrals}
\smallskip

We have:

\begin{align}\label{rel:Hhessian}
\frac{\partial^2 H}{\partial \lambda_i\partial\lambda_j}=\left(-\eps\frac{\partial^2 F^{BM}}{\partial\lambda_i\partial\lambda_j}+\eps^2\frac{\partial F^{BM}}{\partial\lambda_i}\frac{\partial F^{BM}}{\partial\lambda_j}\right)e^{-\eps F^{BM}}.
\end{align}

Using the definition of $F^{BM}$ in \eqref{def:FBM}, Lemma~\ref{lemma:mudiff} and  assuming $\mu_i,i=1,2,3$, to be differentiable functions of $\lambda_i,i=1,2,3$,  we have:

\be\label{rel:eigenBM}
\frac{\partial F^{BM}}{\partial\lambda_j} =\sum_{i=1}^3 \frac{\partial\mu_i}{\partial\lambda_j}(\lambda_i+\frac{1}{3})+\mu_j-\sum_{i=1}^3 \frac{\partial\ln Z}{\partial\mu_i}\frac{\partial\mu_i}{\partial\lambda_j}=\mu_j,\forall j=1,2,3,
\ee where for the last equality we used \eqref{def:musystem}.

On the other hand taking the partial derivative with respect to $\lambda_j, j=1,2,3$ in \eqref{def:musystem} we have:

\be\label{rel:lnZ2derivinv}
\sum_{k=1}^3\frac{\partial^2\ln Z}{\partial\mu_i\partial\mu_k}\frac{\partial{\mu_k}}{\partial\lambda_j}=\delta_{ij}.
\ee

Let us denote by $N=(N_{ij})_{i,j=1,2,3}$ the matrix with components $N_{ij}:=\frac{\partial^2 \ln Z}{\partial\mu_i\partial\mu_j}, i,j=1,2,3$. Then relations (\ref{rel:eigenBM}--\ref{rel:lnZ2derivinv}) show that

$$\left(\frac{\partial \mu_i}{\partial\lambda_j}\right)_{i,j=1,2,3}=\left(\frac{\partial^2 F^{BM}}{\partial\lambda_i\partial\lambda_j}\right)_{i,j=1,2,3}=N^{-1}.$$
 Hence, taking into account \eqref{rel:Hhessian}-\eqref{rel:lnZ2derivinv}, in order to show that $H$ is concave we need to prove that for a suitable $\eps>0$ the matrix:

$$N^{-1}-\eps \mu\otimes \mu$$ is positive definite (where we denoted by $\mu\otimes \mu$ the matrix of components $\mu_i\mu_j, j=1,2,3$).

Let us recall now (see for instance \cite{Bhatia}, Proposition $1.2.6$) that if $P,Q$ are real-valued symmetric matrices with $Q$ and $PQ+QP$ positive definite then $P$ is positive definite as well. We aim to apply this criterion to $Q=N$ and $P=N^{-1}-\varepsilon\mu\otimes \mu$. To this end we need to show that $N$ is positive definite and

\begin{equation} \label{poss}
 \frac{2}{\eps}{\rm Id}-\left(N\mu\otimes\mu+\mu\otimes\mu N\right)\textrm{ is positive definite}
\end{equation} for suitable $\eps>0$. The proof of these two claims is quite technical and  is postponed to the Appendix.

\eproof


\section{The case of a singular heat flux law}
\label{sec:singheat}

In this last section, we consider the case when the heat
flux law exhibits a singular behavior for
$\vt\sim 0$. Namely, in place of \eqref{h12},
we ask that
\begin{equation}\label{h12s}
  \kappa(\vt) = A_0 + A_k \vt^k + A_{-2} \vt^{-2}, \ \
    A_0, A_k, A_{-2} > 0.
\end{equation}
Indeed other types of singular behavior may be considered
as well. We chose the above expression since, as noted in the
Introduction, it is in agreement with the behavior
observed in various types of phase-transition
and phase-separation models both related to
liquids and to solids.

Mathematically speaking, the above choice permits to
consider a more general class of singular potentials
(still including the Ball-Majumdar case). Moreover,
it permits us to get a stronger version of the entropy
inequality (namely, it holds as a distributional inequality,
and not only in the integrated form \eqref{i28weak}).
As a drawback, we have a marginal regularity loss for $\QQ$.
To state our related result, we first have to introduce
the auxiliary function (cf.~\eqref{h12})
\begin{equation}\label{defiH}
  H(\vt) := A_0 \log \vt
   + \frac{A_k}k \vt^k
   - \frac{A_{-2}}2 \vt^{-2},
\end{equation}
needed in the statement of the entropy inequality.

\greybox{
\bete\label{th:main:2}
 Let the coefficients $\mu,\kappa,\Gamma$ satisfy~\eqref{h11},
 \eqref{h13}, and~\eqref{h12s}, with relation~\eqref{defiA}
 on exponents. Moreover,
 let the initial data comply with\/ \eqref{m6}-\eqref{m6c}.
 Finally, let $f$ be any potential fulfilling conditions
 {\rm (i)-(iii)} of Subsection~\ref{subsec:hyp}.
 Then, there exists a quadruplet $(\vu,p,\QQ,\vt)$
 with the regularity\/ \eqref{rego:u}, \eqref{rego:p}
 and\/ \eqref{rego:teta}, together with 
 \begin{align}\label{rego:Qs}
   & \QQ \in L^\infty(0,T;H^1(\Omega;\TS))\cap
    L^{9/5}(0,T;W^{2,9/5}(\Omega;\TS)), \ \
    \QQ_t \in L^{\frac{18}{11}}((0,T)\times\Omega;\TS),\\
  \label{rego:teta-1s}
   & \log \vt \in L^\infty(0,T;L^1(\Omega)) \cap
    L^2(0,T;H^1(\Omega)), \ \
     \vt^{-1} \in L^2(0,T;H^1(\Omega)),
 \end{align}
 satisfying the incompressibility condition~\eqref{i21},
 the momentum equation~\eqref{i22}, the Q-tensor
 evolution system~\eqref{i23} and the total energy
 balance~\eqref{i24} in the sense of distributions
 and complying with the initial conditions as specified
 in~\eqref{init:w0}-\eqref{init:w}.
 Moreover, the entropy inequality holds in the 
 distributional sense: for any nonnegative function
 $\phi \in \calD([0,T]\times\overline{\Omega})$,
 \begin{align}\label{i28we}
  &  \int_0^T\io s \phi_t
    + \int_0^T\io s \vu \cdot \Grad \phi
    +\int_0^T \io H(\vt) \Delta \phi\\
  \no
  & \mbox{}~~~~~
     \le -\int_0^T\io \frac{\phi}{\vt} \left( \frac{\mu(\vt)}{2}
        \left| \Grad \vu + \Grad^t \vu \right|^2
     + \Gamma(\vt) | \HH |^2
     + \frac{\kappa(\vt)}{\vt} |\Grad \vt |^2 \right).
 \end{align}
\ente
}
\beos\label{weaktostrong}
 Let us notice that under the current hypotheses, if
 we have any {\sl sufficiently smooth}\/ weak solution
 satisfying \eqref{i28we} {\sl with the equal sign},
 then that solution also satisfies the heat
 equation \eqref{p6}, or, in other words, it solves the problem
 in the ``physical'' sense. Hence, the current formulation
 and the physical one are in some sense equivalent.
 Instead, the same does not hold, even for
 smooth solutions, when we only know
 that the {\sl integrated}\/ entropy equality
 (i.e., \eqref{i28} with the equal sign) holds.
\eddos
\beos\label{f:admiss}
 Regarding the admissible class of potentials,
 we may relax even more our assumptions. Indeed, what
 we need for $f$ is simply that to be a convex function,
 with open bounded domain $\calD[f]$, smooth inside $\calD[f]$.
 Namely, we do not need $\calD[f]$ to be related to the
 eigenvalues of $\QQ$ (cf.~(ii) of Subsection~\ref{subsec:hyp}).
\eddos

\smallskip

\noindent%
{\bf Proof of Theorem~\ref{th:main:2}.}~~%
The argument follows the lines of the proof given
for Theorem~\ref{th:main}. Hence, we only sketch the
differences occurring in the a-priori estimates and in
the weak stability argument.

First of all, we can observe that all the a priori bounds
directly deriving from the energy and entropy balances
still hold. Additionally, thanks to \eqref{h12s}, \eqref{st15}
is improved as
\begin{equation}\label{st15s}
  \| \Grad \vt^{\kappa/2} \|_{L^2(0,T;L^2(\Omega;R^3))}
   + \| \Grad \log\vt \|_{L^2(0,T;L^2(\Omega;R^3))}
   + \| \Grad \vt^{-1} \|_{L^2(0,T;L^2(\Omega;R^3))}
  \le c.
\end{equation}
Combining the last of \eqref{st15s} with \eqref{st12a} and
using a generalized version of Poincar\'e's inequality (cf. \cite[Lemma~3.2]{GSZ} or \cite[Lemma~5.1]{Kenm}), we
deduce more precisely
\begin{equation}\label{st15s2}
 \| \vt^{-1} \|_{L^2(0,T;H^1(\Omega;R^3))}
  \le c.
\end{equation}
The key point, as before, is constituted by the estimate of
the singular potential.

\smallskip

\noindent%
{\bf Estimate of the singular potential for singular
  heat flux.}~~%
Differently from before, we  do not  need to rely on the explicit
properties of the Ball-Majumdar potential. Indeed, we may directly compute
\begin{align}\label{keyf2}
  & \| \vt^{-1/2} \HH \|_{L^2(\Omega;\TS)}^2
   = \io \vt^{-1} | \Delta \QQ |^2
    + \lambda^2 \io \vt^{-1} | \QQ |^2
    + \io \vt \left| \calL \left[ \frac{\partial f(\QQ)}{\partial \QQ } \right] \right|^2 \\
 \nonumber
  & \mbox{}~~~~~
    + 2\lambda \io \vt^{-1} \Delta \QQ : \QQ
    - 2 \io \calL \left[ \frac{\partial f(\QQ)}{\partial \QQ } \right] :
     \Delta \QQ
     - 2 \lambda \io \calL \left[ \frac{\partial f(\QQ)}{\partial \QQ } \right] : \QQ
    =: \sum_{j=1}^6 T_j.
\end{align}
Now, let us observe that
\begin{align}\label{test21}
  | T_4 | & = 2\lambda \left| \io \vt^{-1} \Delta \QQ : \QQ \right|
   \le \frac12 \io \vt^{-1} | \Delta \QQ |^2
    + 2\lambda^2 \io \vt^{-1} | \QQ |^2, \\
 \label{test22}
   T_5 & =  - 2 \io \calL \left[ \frac{\partial f(\QQ)}{\partial \QQ } \right] :
        \Delta \QQ
      = 2 \io \frac{\partial^2 f(\QQ)}{\partial \QQ_{ij} \partial \QQ_{kl} } \Grad \QQ_{ij}
         \cdot \Grad \QQ_{kl},\\
 \label{test23}	
  | T_6 | & = 2 \lambda \left| \io \calL \left[ \frac{\partial f(\QQ)}{\partial \QQ } \right] : \QQ \right|
   \le \frac12 \io \vt \left| \calL \left[ \frac{\partial f(\QQ)}{\partial \QQ } \right] \right|^2
    + 2\lambda^2 \io \vt^{-1} | \QQ |^2.
\end{align}
Then, integrating \eqref{keyf2} in time over $(0,T)$ and using
the latter of \eqref{st17}, we arrive at
\begin{align}\label{test24}
  & \iTo \vt^{-1} | \Delta \QQ |^2
   + \iTo \frac{\partial^2 f(\QQ)}{\partial \QQ_{ij} \partial \QQ_{kl} } \Grad \QQ_{ij}
         \cdot \Grad \QQ_{k,l}\\
   \nonumber
  & \mbox{}~~~~~
    + \io \vt \left| \calL \left[ \frac{\partial f(\QQ)}{\partial \QQ } \right] \right|^2
   \le c + c \io \vt^{-1} | \QQ |^2
   \le c,
\end{align}
the last inequality following from estimates \eqref{st18} and
\eqref{st15s2}. Hence, in place of \eqref{st21a}-\eqref{st21} we
get now
\begin{equation}\label{st21s}
  \left\| \vt^{1/2} \calL \left[ \frac{\partial f(\QQ)}{\partial \QQ } \right] \right\|_{L^2((0,T)\times \Omega;\TS)}
  + \big\| \vt^{-1/2} \Delta \QQ \big\|_{L^2((0,T)\times \Omega;\TS)}
    \le c.
\end{equation}
Using \eqref{st19} and interpolation, we then obtain
\begin{equation}\label{st21s2}
  \left\| \vt \calL \left[ \frac{\partial f(\QQ)}{\partial \QQ } \right]
    \right\|_{L^{\frac{2A}{A+1}}((0,T)\times \Omega;\TS)}
  + \left\| \Delta \QQ \right\|_{L^{\frac{2A}{A+1}}((0,T)\times \Omega;\TS)}
    \le c,
\end{equation}
where $\frac{2A}{A+1}>\frac95$ thanks to \eqref{defiA}.

\smallskip

\noindent%
{\bf Weak sequential stability for singular heat flux.}~~%
Also this part of the procedure goes through as before, the only
difference being represented by some loss of regularity on $\QQ$.
Indeed, in place of \eqref{lim12}, we now have, due to \eqref{st21s2},
\begin{equation}\label{lim12s}
  \QQ_n \to \QQ \ \ \text{weakly in }
   L^{\frac95+}(0,T;W^{2,\frac95+}(\Omega;\TS)).
\end{equation}
Combining \eqref{st21s} with the second \eqref{lim12b}
and using the Gagliardo-Nirenberg inequalities,
it is not difficult to obtain
\begin{equation}\label{lim12s3}
  \QQ_n \to \QQ \ \ \text{weakly in }
   L^{\frac{18}5+}(0,T;W^{1,\frac{18}5+}(\Omega;\TS)).
\end{equation}
This leads to some obvious modifications
in the exponents in the subsequent relations. For instance, in
\eqref{lim18}-\eqref{lim19}
$L^{12/7}$ turns to $L^{18/11}$; in the
latter of \eqref{lim20b}, $L^2$ has to be replaced,
say, by $L^{9/5}$; in \eqref{lim31}, $L^{4/3}$ becomes
$L^{9/7}$. Finally, $L^{36/29}$
is substituted by $L^{6/5}$ in \eqref{lim35}.
A further check shows that the argument used
to get the strong convergence \eqref{lim43c}
still works even though the summability properties
of the term $|\Grad \QQ|^2/2$ in the expression of
$e$ \eqref{defie} are now a little bit lower.

The subsequent convergence relations regarding $\vt$
do not change. Concerning, instead, the identification
of the singular potential term, it suffices to notice
that, in place of \eqref{lim51}, we now have only
\begin{equation}\label{lim51s}
  \vt_n \calL \left[ \frac{\partial f(\QQ_n)}{\partial \QQ_n } \right]
   \to \ov{\vt \calL \left[ \frac{\partial f(\QQ)}{\partial \QQ } \right]} \ \
  \text{weakly in} \ L^{\frac95}((0,T)\times \Omega;\TS),
\end{equation}
due to \eqref{st21s2}.

To conclude the proof, we have to see that we can take
the (supremum) limit in the distributional entropy
inequality \eqref{i28we} (of course, we are now
assuming it holds in that form at the level $n$
and want to show that it holds in
the same form also in the limit).
To this aim, we first notice that
the \rhs\ is treated as in \eqref{i28lim}, taking now
a generic (nonnegative) test function $\phi$.
To proceed, we have to deal with the \lhs. We notice
first that, by  \eqref{st15s2},
\begin{equation}\label{lim14a}
  \vt_n^{-1} \to \ov{\vt^{-1}} \ \
   \text{weakly in } L^2(0,T;H^1(\Omega)),
\end{equation}
whence, thanks to pointwise convergence (see \eqref{lim45})
and Sobolev's embeddings,
\begin{equation}\label{lim46}
  \vt_n^{-1} \to \vt^{-1} \ \
   \text{strongly in} \ L^p(0,T;L^q(\Omega))
   \ \text{for all }p\in [1,2), \ q\in [1,6).
\end{equation}
In particular, we identify $\ov{\vt^{-1}} = \vt^{-1}$.

Next, let us note that, by \eqref{st21s}, \eqref{lim14a}
and interpolation,
\begin{equation}\label{test31}
  \calL \left[ \frac{\partial f(\QQ_n)}{\partial \QQ_n } \right]
    \to \calL \left[ \frac{\partial f(\QQ)}{\partial \QQ } \right] \ \
   \text{weakly in} \ L^{4/3}(0,T;L^{12/7}(\Omega;\TS)).
\end{equation}
Let us now recall the standard identity
\begin{equation}\label{test32}
  \mathcal{L} \left[ \frac{\partial f (\QQ)} {\partial \QQ} \right] : \QQ
   = f(\QQ) + f^* \left( \mathcal{L} \left[ \frac{\partial f (\QQ)} {\partial\QQ} \right] \right),
\end{equation}
where $f^*$ is the {\sl convex conjugate}\/ (see, e.g., \cite[Thm. 23.5, p. 218]{Rocka}) of $f$.
Then, by (ii) of Subsection~\ref{subsec:hyp},
$f^*$ grows at most linearly at
infinity. Combining \eqref{test31}, \eqref{test32} and \eqref{lim12b},
we then obtain
\begin{equation}\label{test33}
  f(\QQ_n)\to f(\QQ) \ \
   \text{weakly in} \ L^{4/3}(0,T;L^{12/7}(\Omega)).
\end{equation}
Now, interpolating exponents in \eqref{lim11} and accounting
for the strong convergence \eqref{lim23}, we get
\begin{equation}\label{test34}
  \vu_n \to \vu \ \
   \text{strongly in} \ L^{4+}(0,T;L^{\frac83+}(\Omega;R^3)).
\end{equation}
Moreover, we have
\begin{equation}\label{lim37w}
  \| \vt_n^{m-1} \|_{ L^{\frac32}((0,T)\times\Omega)} \le c.
\end{equation}
Then, recalling the expression \eqref{entropy} for $s$,
combining \eqref{test33} with \eqref{test34}, and
\eqref{lim37w} and \eqref{lim45b} with \eqref{lim23},
we finally infer
\begin{align}\label{lim37w2}
  & s_n \to s
   \ \ \text{say, strongly in}\ L^{1+}((0,T)\times\Omega),\\
 \label{lim37w3}
  & s_n \vu_n \to s \vu
   \ \ \text{strongly in}\ L^{1+}((0,T)\times\Omega;R^3).
\end{align}
To conclude the proof, we need to take the limit in the
third integral in the first row of \eqref{i28}, where we
recall that $H$ was defined in \eqref{defiH}.
To this aim,
we need a proper compactness tool:
\bele\label{lemma:comp}
 Let $h:\RR\to [0,+\infty)$ be a continuous function such
 that $h(z)\to\infty$ as $z\to\infty$.
 Suppose that $\{u_n\}$ is a sequence of real-valued
 measurable functions defined over $(0,T)\times \Omega$ such that
 \begin{align}\label{un1}
   & u_n(t,x) \ge 0\ \
    \text{for a.e.}\ x\in (0,T)\times \Omega,\\
  \label{un2}
   & {\rm ess}\sup_{t\in (0,T)}\|h(u_n)(t,\cdot)\|_{L^1(\Omega)}\leq H,\quad
     \|u_n\|_{L^2(0,T;L^p(\Omega))}\leq C\hbox{ for some $p>2$}, \\
  \label{un3}
   & u_n \to u \ \
    \text{a.e. in}\ (0,T)\times \Omega,
 \end{align}
  for some positive constants $H$ and $C$. Then
 \begin{equation}\label{unlim}
   u_n \to u \ \
    \text{strongly in}\ L^2((0,T)\times \Omega).
 \end{equation}
\enle
\begin{proof}
Clearly, it is enough to show equi-integrability of the $L^2$-norm of $u_n$, specifically
\begin{equation}\label{l1}
  \int\int_{u_n\geq M}|u_n|^2\, {\rm dx}\, {\rm dt}\to 0\quad\hbox{as }M\to \infty.
\end{equation}
To see this, we first observe that
\[
  H\geq \int_{u_n(t,\cdot)\geq M} h(u_n)(t, \cdot)\,
    {\rm dx}\geq \inf_{z\geq M}h(z)\left|\{u_n\geq M\}\right|\quad\hbox{for a.e. }t\in (0,T).
\]
Now, applying the H\"older inequality, we get
\[
  \int_{u_n(t,\cdot)\geq M}u_n^2(t,\cdot)\, {\rm dx}
   \leq \left|\{u_n(t, \cdot)\geq M\}\right|^{1/q}\|u_n(t,\cdot)\|_{L^p(\Omega)}^2,
    \quad q=\frac{p}{p-2},
\]
whence \eqref{l1} follows.
\end{proof}
\noindent%
Relations \eqref{st12a} and \eqref{lim14a} allow us to
apply the lemma with the choice $u_n=1/\vt_n$ and $h(u)=\log u$.
Hence, relation \eqref{lim46} is improved up to
\begin{equation}\label{lim38b}
  \vt_n^{-1} \to \vt^{-1} \ \
   \text{strongly in} \ L^2((0,T)\times \Omega).
\end{equation}
Hence, using \eqref{lim45}, we have that
\begin{equation}
  H(\vt_n) \to H(\vt) \ \ \text{strongly in}\ L^{1}((0,T)\times\Omega;R^3).
\end{equation}
and consequently all terms on the left hand side of \eqref{i28we}
pass to the limit $n\nearrow\infty$, which concludes
the proof of the theorem.


\begin{appendix}

\section{A technical estimate}

We conclude the paper by providing the proof of the key relation (\ref{poss}), along with the necessary technical preparations. Our approach is based 
on the co-called Laplace's method to evaluate the integrals containing exponentials, see \cite{second-laplace}.
Let us first note that the positive definiteness of $N=\left(\frac{\partial^2\ln Z}{\partial\mu_i\partial\mu_j}\right)_{i,j=1,2,3}$ is proved in \cite{ball-lecture}. For convenience we recall the argument, namely:

\begin{align*}
\frac{\partial^2\ln Z}{\partial\mu_i\partial\mu_j}a_i a_j&=\frac{1}{Z}\int_{\mathbb{S}^2} e^{\sum_{k}\mu_kp_k^2}a_i p_i^2 a_j p_j^2\,{\rm d}p-\frac{1}{Z^2}\left(\int_{\mathbb{S}^2} e^{\sum_k\mu_kp_k^2}a_i p_i^2\,{\rm d}p\right)\left(\int_{\mathbb{S}^2} e^{\sum_k\mu_kp_k^2}a_j p_j^2\,{\rm d} p\right)\nonumber\\
&=\frac{\int_{\mathbb{S}^2}e^{\sum_k\mu_kp_k^2}a_i p_i^2 a_j p_j^2\,{\rm d}p \int_{\mathbb{S}^2}e^{\sum_k\mu_kq_k^2}{\rm d} q}{\int_{\mathbb{S}^2}\int_{\mathbb{S}^2}e^{\sum_k\mu_k(p_k^2+q_k^2)}\,{\rm d}p\,{\rm d}q}-\frac{\int_{\mathbb{S}^2} e^{\sum_k\mu_kp_k^2}a_i p_i^2\,{\rm d}p\int_{\mathbb{S}^2} e^{\sum_k\mu_kq_k^2}{a_jq_j^2}\,{\rm d} q}{\int_{\mathbb{S}^2}\int_{\mathbb{S}^2}e^{\sum_k\mu_k(p_k^2+q_k^2)}\,{\rm d}p\,{\rm d}q}\nonumber\\
&=\frac{1}{2Z^2}\int_{\mathbb{S}^2}\int_{\mathbb{S}^2} e^{ \sum_k\mu_k(p_k^2+q_k^2)}(a_ip_i^2-a_jq_j^2)^2\,{\rm d}p\,{\rm d}q>0,\, \forall (a_1,a_2,a_3)\not=(0,0,0).
\end{align*}

In estimating the positivity of $\frac{2}{\eps}Id- \left(N\mu\otimes\mu+\mu\otimes\mu N\right)$ we consider the components of the matrix $N\mu\otimes\mu+\mu\otimes\mu N$, all functions of $(\mu_1,\mu_2,\mu_3)\in X$  (cf. \eqref{def:X}):

\begin{align}
\mcI_{ij}&=\frac{\int_{\mathbb{S}^2} e^{\sum_m\mu_m p_m^2} (\sum_k\mu_k p_k^2)(\mu_ip_j^2+\mu_jp_i^2)\,{\rm d}p}{\int_{\mathbb{S}^2} e^{ \sum_m\mu_m p_m^2}\,{\rm d}p}\nonumber\\
\nonumber
&\quad-\left(\frac{\int_{\mathbb{S}^2} e^{ \sum_m\mu_m p_m^2}  (\sum_k\mu_k p_k^2)\,{\rm d}p}{\int_{\mathbb{S}^2} e^{\sum_m\mu_m p_m^2}\,{\rm d}p}\right)\left(\frac{\int_{\mathbb{S}^2} e^{\sum_m\mu_m p_m^2} (\mu_i p_j^2+\mu_j p_i^2)\,{\rm d}p}{\int_{\mathbb{S}^2} e^{\sum_m\mu_m p_m^2}\,{\rm d}p}\right)\nonumber\\
&=\frac{\int_{\mathbb{S}^2}\int_{\mathbb{S}^2} e^{ \sum_m\mu_m(p_m^2+q_m^2)}( \sum_k\mu_kp_k^2-\sum_l\mu_lq_l^2)\left[\mu_i(p_j^2-q_j^2)+\mu_j(p_i^2-q_i^2)\right]\,{\rm d}p{\rm d}q}{2\int_{\mathbb{S}^2}\int_{\mathbb{S}^2} e^{\sum_m\mu_m(p_m^2+q_m^2)}\,{\rm d}p{\rm d}q}.
\end{align}

It will suffice to show that $\mcI_{ij}(\mu_1,\mu_2,\mu_3), i,j=1,2,3$ is bounded independently of $(\mu_1\mu_2,\mu_3)\in X$. To this end we consider two cases: when  $\mu_i\not=\mu_j$ for $i\not=j$ and when two out the three $\mu_i,i=1,2,3$ are equal.

\subsection{Case I: $\mu_i\not=\mu_j$ for $i\not=j$}

We let $\mu_i=\rho\gamma_i$ with $\rho\ge 0$ and $(\gamma_1,\gamma_2,\gamma_3)\in\mathbb{S}^2\cap X$. It will suffice to prove:

\be
\lim_{\rho\to\infty} \mcI_{ij}(\rho\gamma_1,\rho\gamma_2,\rho\gamma_3)<\infty
\ee  uniformly in $(\gamma_1,\gamma_2,\gamma_3)\in\mathbb{S}^2\cap X$.

We first recall the Laplace's asymptotic expansion method � for evaluating integrals (see for instance \cite{second-laplace} and also \cite{asymptotic-integrals}), namely, let us consider the integral:

$$J(h)=\int_\Omega e^{-h f(x)}g(x)\,dx$$ with $\Omega\subset\mathbb{R}^d$ with $\Omega$ an open bounded set and $f,g\in C^\infty(\overline\Omega,\mathbb{R})$ . 

If $f$ has a single strict global minimum point in $\bar x\in\Omega$ such that $f(\bar x)=0$ then for any $n\in\mathbb{N}$ one has:

\be\label{expansion:Laplace}
J(h)=\sum_{k=0}^n c_k h^{-\frac{k+d}{2}}+o(h^{-\frac{n+d}{2}})
\ee as $h\to\infty$ with $c_k$ coefficients that are explicitly computable in terms of $f$ and $g$. More precisely we will only need the first three, $c_0$, $c_1$, $c_2$ given by

\be\label{def:c0}
c_0=\Gamma (\frac{d+1}{2})\int_{\mathbb{S}^{d-1}} \frac{g_0(\sigma)}{d(f_0(\sigma))^{\frac{d}{2}}}\,{\rm d}\sigma,
\ee

\be\label{def:c1}
c_1=\Gamma(\frac{d+2}{2}) \int_{\mathbb{S}^{d-1}} \frac{2d g_1(\sigma)f_0(\sigma)-(d+1)g_0(\sigma)f_1(\sigma)}{2d(d+1)f_0(\sigma)^{\frac{d+3}{2}}}\,{\rm d}\sigma,
\ee

\be\label{def:c2}
c_2=\Gamma(\frac{d+3}{2})\int_{\mathbb{S}^{d-1}}\frac{1}{d+2}\frac{g_2}{f_0^{\frac{d+1}{2}}}-\frac{f_1g_1}{(d+1)f_0^{\frac{d+4}{2}}}+\frac{f_1^2g_0}{4df_0^{\frac{d+6}{2}}}+\frac{(5f_1^2-4f_0f_2)g_0}{8df_0^{\frac{d+6}{2}}}\,{\rm d}\sigma,
\ee where we denoted, for $\sigma=\frac{x}{|x|}$:

$$f_0(\sigma):=\frac{1}{2}\left(\frac{\partial^2 f}{\partial x_i\partial x_j}(\bar x)\right)\left(\frac{(x-\bar x)_i}{|x-\bar x|}\frac{(x-\bar x)_j}{|x-\bar x|}\right),$$

$$f_1(\sigma)=\frac{1}{6}\left(\frac{\partial^3 f}{\partial x_i\partial x_j\partial x_k}(\bar x)\right)\left(\frac{(x-\bar x)_i}{|x-\bar x|}\frac{(x-\bar x)_j}{|x-\bar x|}\frac{(x-\bar x)_k}{|x-\bar x|}\right),$$

$$g_0(\sigma)=g_0(\bar x), g_1(\sigma)=\frac{\partial g}{\partial x_i}(\bar x)\frac{(x-\bar x)_i}{|x-\bar x|},g_2(\sigma)=\frac{1}{2}\frac{\partial^2 g}{\partial x_i\partial x_j}(\bar x)\frac{(x-\bar x)_i}{|x-\bar x|}\frac{(x-\bar x)_j}{|x-\bar x|}.$$

Let us observe that in \cite{second-laplace} the formulae were determined by allowing $g$ to be potentially singular at the maximum point. However, if $g$ is smooth (as is our case) then one can see that the odd-index coefficients $c_{2k+1}$ are zero (see also \cite{asymptotic-integrals} or directly compute $c_1$ by the formula above).

In order to apply the previous argument let us observe that the function $\gamma_m (p_m^2+q_m^2)$ of variable $(p,q):=((p_1,p_2,p_3),(q_1,q_2,q_3))\in\mathbb{S}^2\times\mathbb{S}^2$  attains its maximum value at two points, that depend on the maximum element in the set $\{\gamma_1,\gamma_2,\gamma_3\}$. Let us assume without loss of generality that $\gamma_2,\gamma_3< \gamma_1$. Then the maximum of the function $\gamma_m (p_m^2+q_m^2)$ is attained at two points $(p,p)$ with $p\in \{(1,0,0),(-1,0,0)\}$ so in order to apply  Laplace's method we need to split $\mathbb{S}^2\times \mathbb{S}^2$ into two subdomains. Let us denote $\mathbb{S}_E:=\{p\in\mathbb{S}^2, p\cdot (1,0,0)>0\}$. We then apply the previously mentioned Laplace's method on each of the sets $\mathcal{E}:=\mathbb{S}_E\times\mathbb{S}_E$ and $\mathcal{V}:=\mathbb{S}^2\times\mathbb{S}^2\setminus\mathcal{E}$ chosing $h=\rho$ and $f(p_1,p_2,p_3,q_1,q_2,q_3):=\gamma_1-\gamma_i (p_i^2+q_i^2)$ (note that we can multiply both denominator and numerator in $\mcI_{ij}$ by $e^{-\rho\gamma_1}$).

The function $g$ will be chosen to be either $g\equiv 1$ (for treating denominators) or
\be\label{def:gnumerator}
g_{ij}(p_1,p_2,p_3,q_1,q_2,q_3):=(\gamma_kp_k^2-\gamma_lq_l^2)\left[\gamma_i(p_j^2-q_j^2)+\gamma_j(p_i^2-q_i^2)\right]
\ee
(for dealing with the numerators). Then we take $d=4$ and  for treating the denominator we note that each of the two integrals  (over $\mathcal{E}$, respectively $\mathcal{V}$)  admit an asymptotic expansion of the form: $\bar c_0\rho^{-2}+o(\rho^{-2})$. Similarily for the numerator the integrals  have an expansion of the form  $\rho^2(\tilde c_0\rho^{-2}+\tilde c_1\rho^{-\frac{5}{2}}+\tilde c_2\rho^{-3}+\tilde c_3\rho^{-\frac{7}{2}}+\tilde c_4\rho^{-4}+o(\rho^{-4}))$. As previously mentioned the coefficients of odd index, namely $\tilde c_1,\tilde c_3$ are zero. Taking into account the specific form of $g_{ij}$ \eqref{def:gnumerator} and the expressions of the coefficients $\tilde c_0$, namely \eqref{def:c0} respectively $\tilde c_2$, namely \eqref{def:c2} we have that both $g_{ij}$ and its first and second derivatives are zero\footnote{note that the derivatives are computed at $(\theta,\varphi,\tilde\theta,\tilde\varphi)=(0,\frac{\pi}{2},0,\frac{\pi}{2})$ and we take  $(p_1,p_2,p_3)=(\cos\varphi\sin\theta,\sin\varphi\sin\theta,\cos\theta)$ with $(\varphi,\theta)\in [0,2\pi]\times (0,\pi)$ and $(q_1,q_2,q_3)=(\cos\tilde\varphi\sin\tilde\theta,\sin\tilde\varphi\sin\tilde\theta,\cos\tilde\theta)$ with $(\tilde\varphi,\tilde\theta)\in [0,2\pi]\times (0,\pi)$}  so $\tilde c_0=\tilde c_2=0$. Thus the numerator will have an asymptotic expansion of the form $\tilde c_4\rho^{-2}+o(\rho^{-2})$ and overall $\mcI_{ij}$ will be bounded as $\rho\to\infty$.

\bigskip
\subsection{Case II:  two out the three $\mu_i,i=1,2,3$ are equal}

We assume without loss of generality that $\mu_1=\mu_2$ and we denote the common value by $\mu$. Then $\mu_3=-\mu_1-\mu_2=-2\mu$. After a couple of manipulations, using that $(p_1,p_2,p_3)=(\cos\varphi\sin\theta,\sin\varphi\sin\theta,\cos\theta)$ with $(\varphi,\theta)\in [0,2\pi]\times (0,\pi)$ we have (by denoting $\bZ:=2\pi\int_0^\pi \emu \sin\theta\,d\varphi$):

\begin{align*}
\mcI_{11}=\frac{6\mu^2}{\bZ}\int_0^{2\pi}\int_0^\pi \emu &\cos^2\varphi\sin^5\theta\,{\rm d}\varphi {\rm d}\theta\\
&-\frac{12\pi\mu^2}{\bZ^2}\int_0^\pi \emu\sin^3\theta\,{\rm d}\theta \int_0^{2\pi}\int_0^\pi \emu\cos^2\varphi \sin^3\theta\,{\rm d}\varphi {\rm d}\theta,
\end{align*}

\begin{align*}
\mcI_{22}=\frac{6\mu^2}{\bZ}\int_0^{2\pi}\int_0^\pi \emu &\sin^2\varphi\sin^5\theta\,{\rm d}\varphi {\rm d}\theta\\
&-\frac{12\pi\mu^2}{\bZ^2}\int_0^\pi\emu \sin^3\theta\,d\theta \int_0^{2\pi}\int_0^\pi \emu \sin^2\varphi \sin^3\theta\,{\rm d}\varphi {\rm d}\theta,\
\end{align*}

\be
\mcI_{33}=4\mcI_{12}=4\mcI_{21}=12\mu^2\left[\frac{2\pi}{\bZ}\int_0^\pi \emu \sin^5\theta\,{\rm d}\theta-\left(\frac{2\pi\int_0^\pi \emu\sin^3\theta\,{\rm d}\theta}{\bar Z}\right)^2\right],
\nonumber
\ee

\begin{align*}
\mcI_{13}=\mcI_{31}=\frac{3\mu^2}{\bZ}\int_0^{2\pi}\int_0^\pi &\emu\sin^5\theta(-1-2\cos^2\varphi)\,{\rm d}\theta\,{\rm d}\varphi\\
&-\frac{6\pi\mu^2}{\bZ^2}\int_0^\pi \emu \sin^3\theta\,{\rm d}\theta\int_0^{2\pi} \int_0^\pi \emu (-1-2\cos^2\varphi)\sin^3\theta\,{\rm d}\theta\,{\rm d}\varphi,
\end{align*}

\begin{align*}
\mcI_{23}=\mcI_{32}=\frac{3\mu^2}{\bZ}\int_0^{2\pi}\int_0^\pi &\emu\sin^5\theta(-1-2\sin^2\varphi)\,{\rm d}\theta\,{\rm d}\varphi\\
&-\frac{6\pi\mu^2}{\bZ^2}\int_0^\pi \emu \sin^3\theta\,d\theta \int_0^{2\pi}\int_0^\pi \emu (-1-2\sin^2\varphi)\sin^3\theta\,{\rm d}\theta\,{\rm d}\varphi.
\end{align*}

Thus it suffices to understand what happens with ratios of the type $\frac{\int_0^{\pi} \emu (\mu\sin^2\theta)^k \sin\theta\,d\theta}{\int_0^{\pi} \emu\sin\theta\,d\theta}$ for $k=1,2$. The case when $\mu\to\infty$ can be dealt with using Laplace's method previously described (but in $1D$ now, as opposed to $2D$ domains before),  as in this case the maximum of $3\sin^2\theta$ occurs in the interior of the interval $(0,\pi)$.  We obtain as desired that the limit as $\mu\to\infty$ of all the $\mcI_{ij}$ is finite.

In the case $\mu\to -\infty$ we note that $-3\sin^2\theta$ attains it maximum at the endpoints of the interval $(0,\pi)$ so Laplace's method cannot be applied directly. In this case we let $-\mu=\alpha^2$ and consider the change of variables:
$\alpha\sin\theta=y$ which leads to:

$$\frac{\int_0^{\pi}\emu(\mu\sin^2\theta)^k \sin\theta\,{\rm d}\theta}{\int_0^{\pi} \emu\sin\theta\,{\rm d}\theta}=\frac{\int_0^{\frac{\pi}{2}}\emu(\mu\sin^2\theta)^k \sin\theta\,{\rm d}\theta}{\int_0^{\frac{\pi}{2}} \emu\sin\theta\,{\rm d}\theta}=(-1)^k\frac{\int_0^\alpha e^{-3y^2}y^{2k+1} \frac{{\rm d}y}{\sqrt{1-\frac{y^2}{\alpha^2}}}}{\int_0^\alpha e^{-3y^2}{ y}\frac{{\rm d}y}{\sqrt{1-\frac{y^2}{\alpha^2}}}}.
$$

We let $f(\alpha):=\frac{\int_0^\alpha e^{-3y^2}y^{2k+1} \frac{{\rm d}y}{\sqrt{1-\frac{y^2}{\alpha^2}}}}{\int_0^\alpha e^{-3y^2}{y}\frac{{\rm d}y}{\sqrt{1-\frac{y^2}{\alpha^2}}}}$ and then we have:

\begin{equation}\label{eq:Iest0}
f(\alpha)\le\frac{\int_0^{\frac{\alpha}{2}} \left(e^{-3y^2}\right)\frac{y^{2k+1}\,{\rm d}y}{\sqrt{1-\frac{1}{4}}}+e^{-\frac{3\alpha^2}{4}}\alpha^{2k+1}\int_{\frac{\alpha}{2}}^\alpha \frac{dy}{\sqrt{1-\frac{y^2}{\alpha^2}}}}{\int_0^{\frac{\alpha}{2}}e^{-3y^2}{ y}\,{\rm d}y}.
\end{equation}

On the other hand:

\begin{align*}&\int_{\frac{\alpha}{2}}^{\alpha} \frac{{\rm d}y}{\sqrt{1-\frac{y^2}{\alpha^2}}}=\int_{\frac{\alpha}{2}}^{\alpha} \frac{\alpha {\rm d}y}{\sqrt{(\alpha-y)(\alpha+y)}}\le \frac{\alpha}{\sqrt{\frac{3}{2}\alpha}}\int_{\frac{\alpha}{2}}^{\alpha}\frac{{\rm d}y}{\sqrt{\alpha-y}}\\
&\le \sqrt{\frac{2}{3}}\alpha^{\frac{1}{2}}\int_0^{\frac{\alpha}{2}}\frac{{\rm d}z}{\sqrt{z}}=\frac{2}{\sqrt{3}}\alpha.
\end{align*}

We assume without loss of generality that $\alpha>1$\footnote{that $f(\alpha)$ is bounded when $\alpha\in [0,1]$ can be more easily seen by looking at the expression of $f(\alpha)$ in trigonometric form} and using this assumption as well as the last estimate in \eqref{eq:Iest0} we obtain:

\begin{equation}\nonumber
f(\alpha)\le\frac{\int_0^\infty \left(e^{-3y^2}\right)\frac{y^{2k+1}\,{\rm d}y}{\sqrt{1-\frac{1}{4}}}+e^{-\frac{3\alpha^2}{4}}\alpha^{2k+1}\frac{2}{\sqrt{3}}\alpha}{\int_0^{1/2} e^{-y^2}{y}\,{\rm d}y},
\end{equation} so $f(\alpha)$ is bounded independently of $\alpha$.

\end{appendix}


\end{document}